\documentclass[preprint,authoryear, 12pt]{elsarticle}
\usepackage{amsmath,amsthm,amssymb,amsfonts,fullpage}
\usepackage{color}
\usepackage{enumerate}

\newcommand{\sfr}{\color{black}}

\def \bs {\boldsymbol}
\def \mc {\mathcal}

\def \er {P^-}
\def \dt {\mathrm{mDT}}
\def\prob{\mathbb{P}}

\def\expt{\mathbb{E}}
\def\real{\mathbb{R}}

\newcommand{\until}[1]{\{1,\dots, #1\}}

\newcommand{\subscr}[2]{#1_{\textup{#2}}}
\newcommand{\supscr}[2]{#1^{\textup{#2}}}

\newcommand{\seqdef}[2]{\{#1\}_{#2}}

\usepackage{graphicx,subfigure}

\graphicspath{{./SimFigs/}}

\title{A martingale analysis of first passage times of time-dependent Wiener diffusion models}

\author[mae,msu]{Vaibhav Srivastava}
\ead{vaibhav@egr.msu.edu}

\author[ku]{Samuel F. Feng*}
\ead{samuel.feng@kustar.ac.ae}

\author[psyc,pni]{Jonathan D. Cohen}
\ead{jdc@princeton.edu}

\author[mae]{Naomi~Ehrich~Leonard}
\ead{naomi@princeton.edu}

\author[pni]{Amitai~Shenhav}
\ead{amitai_shenhav@brown.edu}

\address[mae]{Department of Mechanical and Aerospace Engineering, Princeton University, Princeton, NJ, USA}

\address[ku]{Department of Applied Mathematics and Sciences, Khalifa University, Abu Dhabi, UAE}

\address[psyc]{Department of Psychology, Princeton University, Princeton, NJ, USA}

\address[pni]{Princeton Neuroscience Institute, Princeton University, Princeton, NJ, USA}

\address[msu]{Department of Electrical and Computer Engineering, Michigan State University, East Lansing, MI, USA}

\address[brown]{Department of Cognitive, Linguistic, \& Psychological Sciences and Brown Institute for Brain Science, Brown University, Providence, RI, USA}

\cortext[cor]{Corresponding author}

\begin{document}
\begin{abstract}
  {\sfr Research in psychology and neuroscience has successfully
    modeled decision making as a process of noisy evidence
    accumulation to a decision bound.  While there are several
    variants and implementations of this idea, the majority of these
    models make use of a noisy accumulation between two absorbing
    boundaries.  } A common assumption of these models is that
  decision parameters, e.g., the rate of accumulation (drift rate),
  remain fixed over the course of a decision, {\sfr allowing the
    derivation of analytic formulas} for the probabilities of hitting
  the upper or lower decision threshold, and the mean decision time.
  There is reason to believe, however, that many types of behavior
  would be better described by a {\sfr model} in which the parameters
  were allowed to vary over the course of the decision process.

In this paper, we use martingale theory to derive formulas for the
mean decision time, hitting probabilities, and first passage time
(FPT) densities of a Wiener process with time-varying drift between
two time-varying absorbing boundaries.  {\sfr This model was first
  studied by \cite{Ratcliff1980} in the two-stage form, and here we
  consider the same model for an arbitrary number of stages
  {(i.e. intervals of time during which parameters are constant).  Our
    calculations enable direct computation of mean decision times and
    hitting probabilities for the associated multistage process.  We
    also provide a review of how martingale theory may be used to
    analyze similar models employing Wiener processes by re-deriving
    some classical results.  In concert with a variety of numerical
    tools already available, the current derivations should encourage
    mathematical analysis of more complex models of decision making
    with time-varying evidence.}}
\end{abstract}

\maketitle

\section{Introduction}
Continuous time stochastic processes modeling a particle's diffusion (with drift) towards one of two absorbing boundaries have been used in a wide variety of applications including statistical physics \citep{Farkas2001}, finance \citep{Lin1998}, economics~\citep{Webb2015}, and health science \citep{Horrocks2004}. 
Varieties of such models  have also been applied extensively within psychology and neuroscience to describe both the behavior and neural activity associated with decision processes involved in perception, memory, attention, and cognitive control~\citep{Heath1992,Ratcliff1998,Ratcliff2008,Simen2009, JIG-MNS:01, JIG-MNS:07, BWB-MMB-CDB:13, SF-PH-AR-WTN:09, MNS-WTN:01, Diederich2014, Diederich2016}; for reviews see \citep{Ratcliff2004, Bogacz2006, Busemeyer2010}. 

In {\sfr these stochastic accumulation decision} models, the state variable $x(t)$ is thought to represent the amount of accumulated noisy evidence at time $t$ for decisions represented by the two absorbing boundaries, that we refer to as the upper ($+$) and lower ($-$) thresholds (boundaries).
The evidence $x(t)$ evolves in time according to a biased random walk with Gaussian increments, which may be written as {\sfr $d x(t) \sim \text{Normal}(\mu \;dt,\sigma ^2 \; dt)$}, and a decision is said to be made at the random time $\tau$, the smallest time $t$ for which $x(t)$ hits either the upper threshold ($x(\tau) = +\zeta$) or the lower threshold ($x(\tau) = -\zeta$), also known as the first passage time (FPT).
The resulting decision dynamics are thus described by the FPT of the underlying {\sfr model}. In studying these processes one is often interested in relating the mean decision time and the probability of hitting a certain threshold (e.g. the probability of making a certain decision) to empirical data.
For example, these metrics can offer valuable insight into how actions and cognitive processes might maximize reward rate, which is a simple function of the FPT properties~\citep{Bogacz2006}.

However, not all decisions can be properly modeled if parameters are fixed throughout the duration of the decision process. 
Certain contexts
can be better described by a {\sfr model} whose parameters {\sfr change with} time.
{\sfr In this article we analyze the time-dependent version of the Wiener process with drift between two absorbing boundaries, building on recent work that is focused on similar time-varying random walk models \citep{Hubner2010,Diederich2014}.
  After reviewing how martingale theory can be used to analyze and re-derive the classical FPT results for the time independent case, we calculate results for the time-dependent case.
  The main theoretical results are presented in \S \ref{subsec:multistage-ddm}, where we provide closed form expressions for threshold-hitting probabilities and expected decision times.
In~\ref{app:ou},
 we also describe how our methods can be applied to the more general Ornstein-Uhlenbeck (O-U) processes, which are similar to the Wiener diffusion processes albeit with an additional ``leak'' term.
We conclude with a summary of the results and a discussion of how the present work interfaces with other similar analyses of time-varying random walk models.}

\section{{\sfr Notation and terminology}}

{\sfr Here we introduce the notation and terminology for describing
  the model we analyze, which is a Wiener process with
  (time-dependent) piecewise constant parameters. This simple
  stochastic model, and others close to it, have been studied before
  \citep{Ratcliff1980, Heath1992, Smith2000, Diederich2003,
    Diederich2006, Diederich2014, Bogacz2006, Wagenmakers2007},
  although the reader should note that our parameterization differs
  from that of some previous studies. Before describing the model, we first review the case where parameters are
  unchanging with time.  In order to easily discuss this simpler model
  alongside the main time-dependent model analyzed in \S
  \ref{sec:performance-metrics}, we use the terms \emph{single-stage
    model} and \emph{multistage model}, respectively.  Readers
  familiar with the popular Diffusion Decision Model of
  \cite{Ratcliff2008} should be aware that parameters in our model do
  not vary randomly trial-to-trial.  Readers familiar with
  \cite{Bogacz2006} should be aware that the single-stage model
  \eqref{eq:pure-ddm} is equivalent to what \cite{Bogacz2006} call the
  ``pure'' Drift Diffusion Model.}

\subsection{The single-stage model/process with constant parameters}\label{subsec:ddm}

{\sfr
  Consider the stochastic differential equation (SDE):
\begin{equation} \label{eq:pure-ddm}
dx(t) = \mu \, dt + \sigma dW(t),\,\,\, x(t_0) = x_0, 
\end{equation}
where parameters $\mu$ and $\sigma$ are constants referred to as the drift and diffusion rates, respectively; $x_0$ is the initial condition (the initial evidence or starting point) of the decision process; and $ \sigma dW(t)$ are independent Wiener increments with variance $\sigma^2 d t$.  This simple stochastic model has successfully modeled the evolution of evidence between two decisions during a two-alternative forced choice task \citep{Ratcliff1998,Bogacz2006, ratcliff2016diffusion}, so that \eqref{eq:pure-ddm} can be interpreted as modeling a decision process in which an agent is integrating noisy evidence until sufficient evidence is gathered in favor of one of the two alternatives. }

A decision is made when the evidence $x(t)$ crosses one of the two
symmetric decision thresholds $\pm \zeta$ for the first time, also
referred to as its \emph{first passage time} (FPT). In other words, a
decision occurs the instance $x(t)$ crosses and is absorbed by one of
the two boundaries. The two boundaries each correspond to one of the
two possible decisions for the task.  We will refer to the absorbing
thresholds at $+\zeta$ and $-\zeta$ as the \emph{upper} and
\emph{lower} decision boundaries/thresholds.  {\sfr To contrast with
  the next section, we will sometimes refer to this model with
  time-invariant parameters as the \emph{single-stage model} or
  \emph{single-stage process}. We again note that the parameterization
  used here differs from that employed by others
  \citep{Smith2000,Ratcliff2004,Navarro2009}, although the underlying
  model is equivalent.}  Our formulation, compared to some others,
does not include a parameter for ``non-decision time'' or ``timeout''
for a given trial. Such a term could be incorporated by shifting the
entire reaction time distribution -- it has no effect on any of our
analyses.




\subsection{Time-dependent, piecewise constant parameters}
\label{sec:msdm}
{\sfr The assumption} that model parameters remain constant throughout the decision process is unlikely to hold in many situations. For example, the quality of evidence may not be stationary (i.e., the drift rate and diffusion rate are not a constant with respect to time) or decision urgency may require thresholds to collapse in order to force a decision by some deadline~\citep{PC-GAP-SE:09, mormann2010drift, zhang2014time, JD-RMB-etal:12}.

In order to analyze such situations, we focus our present study on a two-stage model originally analyzed by \cite{Ratcliff1980}, which we generalize to an arbitrary number of stages. {\sfr We refer to this slightly generalized model as a \emph{multistage model} or \emph{multistage process}, to distinguish it from \eqref{eq:pure-ddm} above.  The multistage model allows for the drift rate, diffusion rate, and thresholds to be piecewise constant functions of time.}

To fully describe the {\sfr multistage model}, we first partition the set of non-negative real numbers (i.e. time axis) into $n$ intervals {\sfr (or \emph{stages})} $[t_{i-1}, t_i], \; i\in \until{n}$ with $t_0=0$ and $t_n = +\infty$. We then assume that the drift rate, the diffusion rate, and the decision thresholds are constant within each interval, but that their values change between intervals.  Evidence {\sfr integration is} thus modeled by
\begin{equation} \label{eq:multistage-ddm}
dx(t) = \mu(t) \, dt + \sigma (t)dW(t),\,\,\, x(t_0) = x_0,
\end{equation}
where 
\begin{align*}
\mu(t) &= \mu_i,  \quad  \text{ for }  t_{i-1} \leq t < t_i, \\
\sigma (t) &= \sigma_i,  \quad   \text{ for } t_{i-1} \leq t < t_i,
\end{align*}
for each $i \in \until{n}$.  The above assumptions are identical to
the assumptions in~\cite{Diederich2014,Diederich2016}.  {\sfr If
  $n=1$, the multistage model reduces to the single-stage
  model~\eqref{eq:pure-ddm}.}  For expository clarity, we begin by
assuming the decision thresholds are fixed at $\pm \zeta$, {\sfr and
  in \S\ref{sec:thresholds} we generalize} to time-varying (piecewise
constant) thresholds. Let $\tau$ be the first passage (decision) time
for the multistage model.

We will frequently refer to the $i$-th stage of \eqref{eq:multistage-ddm}, which for $t>t_{i-1}$, is written as
\begin{equation}\label{eq:ddm-i}
d x(t) = \mu_i dt + \sigma_i dW (t), \,\,\, x(t_{i-1}) = X_{i-1},\,\,\,\text{decision thresholds } \pm \zeta,
\end{equation}
where the initial condition $X_{i-1}$ is a random variable defined as $x(t_{i-1})$ conditioned on there being no decision (threshold-crossing) before time $t_{i-1}$.
More precisely, the density of $X_{i-1}$ is the conditional distribution of $x(t_{i-1})$ given that $\tau > t_{i-1}$. Thus, the random variable $X_{i-1}$ corresponds to realizations of the {\sfr multistage model} that remain within the thresholds $\pm \zeta$ until time $t_{i-1}$. {For the first stage, $X_0$ could be a point mass centered at $x_0$, or it may be a random variable capturing the variability in starting points~\citep{Ratcliff1998}.
The key difference between \eqref{eq:ddm-i} and the single-stage process \eqref{eq:pure-ddm} is that the initial condition $X_{i-1}$ {\sfr is a} random variable whose distribution is determined by previous stages. }

\section{Martingale theory applied to the {\sfr single-stage model}}
\label{sec:review}

  In this section, we give an introduction to the basic properties of
  martingales and the optional stopping theorem, which are the key
  mathematical tools used in calculating our {\sfr main results in \S
    \ref{sec:performance-metrics}.  For readers who are less familiar
    with martingale methods, we first } derive the mean decision time,
  hitting probabilities, and FPT densities for the single- and
  two-stage models. These analyses provide an alternate approach to
  deriving these classical results as compared to {\sfr other}
  non-martingale based approaches
  \citep{Ratcliff1980,Diederich2014}. We discuss the {\sfr differences
    between these various} approaches in \S\ref{sec:discussion}.
  
\subsection{Continuous time martingales}\label{subsec:martingales}
Consider a continuous time stochastic process $\eta(t), \; t  > 0$. Let $\eta(t_1:t_2)$ denote the portion of $\eta(t)$ between times $t_1$ and $t_2$. A stochastic process $M(t)$ is said to be a martingale with respect to $\eta(t)$ if the following three conditions hold:  
\begin{enumerate}
\item $M(t)$ is a function of $\eta(0:t)$ and does not depend on future values of $\eta(t)$\footnote{{\sfr More precisely, $M(t)$ is progressively measurable with respect to the sigma algebra generated by $\eta(0:t)$. See \cite[chap. 2]{Doob1953}, \cite[chap. 1]{Karatzas1998}, or \cite[chap. 1]{Revuz1999}.} }
\item $\expt | M(t)|  < \infty $ 
\item $M(t)$ satisfies {the stationarity condition in expected value} 
  \begin{equation}\label{eq:martingale-definition}
\expt[M(t_2) | \eta(0:t_1)] = M(t_1), \quad \text{for all } 0 \leq t_1 <t_2< \infty.   
\end{equation}
\end{enumerate}

The first condition means that given the realized values $\eta(0:t)$, we should be able to compute $M(t)$ deterministically, so that $M(t)$ does not depend on the future. 
{The second condition is a regularity condition that ensures that $M(t)$ is sufficiently\footnote{$M(t)$ may be heavy tailed due to non-existence of second and higher moments.} light-tailed and holds under several decision-making scenarios.} The third condition is the most crucial -- it enforces stationarity in expected value. This last condition can be interpreted as a ``fair play" condition ensuring that chances of gaining and losing starting with value $M(t_1)$ at time $t_1$ are the same, as in the classic example of a sequence of flips of a fair coin.  When introducing a martingale, one often does not explicitly specify the process $\eta(t)$, and in this case $\eta(t)$ is assumed equal to $M(t)$.

{\sfr Martingale theory is very broad and there are many different choices for $M(t)$ and $\eta$ which are interesting. All of the calculations and results of this paper are constructed from two fundamental stochastic processes, which we now introduce.}
Let $W(t)$ be the standard Wiener process and $X(t) = \mu t + \sigma W(t)$ be a single-stage Wiener process with drift rate $\mu$ and diffusion rate $\sigma$ (without boundaries). We now consider some martingales associated with these two stochastic processes:  
\begin{enumerate}
\item $W(t)$ is a martingale. It is easy to verify that {for $0\le t_1<t_2$},  $\expt[W(t_2)|W(t_1)] = W(t_1)$. Similarly, $X(t) - \mu t$ is a martingale. 
\item $W(t)^2 - t$ is a martingale. Note that {for $0\le t_1<t_2$},  conditioned on $W(t_1)$, $W(t_2)$ has a Gaussian distribution with mean $W(t_1)$ and variance $(t_2-t_1)$. Therefore,
\[
\expt[W(t_2)^2 |W(t_1)] = W(t_1)^2 +(t_2-t_1) \implies \expt[W(t_2)^2 -  t_2 | W(t_1)] = W(t_1)^2 - t_1.
\] 
\item For any $\lambda \in \real$, $\exp(\lambda X(t) -\lambda \mu t - \lambda^2 \sigma^2 t /2)$ is a martingale. Note that {for $0\le t_1<t_2$},  conditioned on $X(t_1)$, $X(t_2)$  has a Gaussian distribution with mean $X(t_1)+\mu (t_2-t_1)$ and variance $\sigma^2(t_2-t_1)$. Thus, 
\begin{align*}
\expt[\exp(\lambda X(t_2))|X(t_1)] &= \exp \left(\lambda (X(t_1)+\mu(t_2-t_1)) + \lambda^2 \sigma^2 (t_2-t_1)/2 \right)\\
\expt[\exp(\lambda X(t_2) - \lambda \mu t_2 -\lambda^2 \sigma^2 t_2/2 )|X(t_1)] &=  \exp \left(\lambda X(t_1) - \lambda \mu t_1 -\lambda^2 \sigma^2 t_1/2 \right).
\end{align*}
For $\lambda = - 2 \mu/\sigma^2$, this martingale reduces to $\exp(-2 \mu X(t)/\sigma^2)$ which is referred to as the \emph{exponential martingale}.
\end{enumerate} 

\subsection{Stopping times and the optional sampling theorem}
\label{subsec:optional}
The first passage time $\tau$ is a random variable defined by $\tau = \inf \{ t>0 | x(t) \notin (-\zeta,\zeta)\}$. 
We are interested in computing conditional expectations and probability densities of $\tau$, which correspond to expected decision times and the corresponding distributions of response times.
The key tool we borrow from the theory of martingales is a classic result known as Doob's optional sampling theorem (also known as the optional stopping theorem), which we motivate and introduce here.
To understand the optional sampling theorem, one must first recall that the expected value of a martingale $M(t)$ computed over all realizations starting from $M(0)$ is equal to the initial expectation of $M(0)$. That is, martingales by definition must satisfy the following:
\[
\expt[M(t)] = \expt[M(0)], \; \text{for each } t >0. 
\]
One then wonders: Does a similar property extend to the random time $\tau$? More specifically, if we consider different realizations of $\tau$ and compute averages of $M(\tau)$ at these realized values, does this average, as the number of realizations grow large, converge to $M(0)$? The answer is affirmative if $\tau$ is well behaved and is independent of the process $M(t)$. Indeed, in this case 
\begin{equation}\label{eq:stopping-motivation}
\expt[M(\tau)|M(0)] = \expt[\expt[M(t)|M(0), \tau=t]] = \expt[M(0)],
\end{equation}
where the outer expectation is with respect to $\tau$ and the inner expectation is with respect to $M(t)$.

But what if $\tau$ is not independent of $M(t)$? In these cases the situation is more subtle.
Suppose $\tau$ is bounded from above by $\bar t$. Then, we can write
{\sfr 
\[
M(\tau) = M(0) + \int_{u=0}^\tau d M(u) = M(0) + \int_{u=0}^{\bar t} \bs 1(\tau \ge u) d M(u),
\]
}
where $\bs 1(\cdot)$ is the indicator function that takes value $0$ if its argument is false and $1$ otherwise, and $dM(u)$ is the (random) increment in $M(u)$ at time $u$. If we assume that the value $\bs 1(\tau \ge u)$ can be deterministically computed based on the knowledge of $M(0:u)$, we can then write
\begin{align*}
\expt[M(\tau)] &=\expt[M(0)] + \int_{u=0}^{\bar t} \expt[\bs 1(\tau \ge u) d M(u)] = 
\expt[M(0)] + \int_{u=0}^{\bar t} \expt[ \expt[\bs 1(\tau \ge s) d M(u) | M(0:u)]], 
\end{align*}
where the second equality follows from the law of total expectation\footnote{For an integrable random variables $Y$ and an arbitrary random variable $Z$, $E[Y] = E[E[Y|Z]]$. Loosely speaking, the law of total expectation states that the expectation of a random variable can be computed by first computing the expectation conditional on another random variable, and then computing the expected value of the resulting expectation. }.  The first equality requires swapping of integral and expectation operators, which is allowed because $\bar t$ is finite. Furthermore, $\bs 1(\tau \ge u) $ is a deterministic function of $M(0:u)$ and thus, $\expt[\bs 1(\tau \ge u) d M(s) | M(0:u)]=\expt[d M(u) | M(0:u)]\bs 1(\tau \ge u) =0$, where the last equality follows by definition of martingale. Consequently, for a random variable $\tau$ and martingale $M(t)$, $\expt[M(\tau)] =\expt[M(0)] $ if (i) the event $\tau\ge u$ is determined by $M(0: u)$, and (ii) $\tau$ is bounded from above with probability one. A random variable satisfying the first condition is called a \emph{stopping time}, and the above discussion is the optional sampling theorem which we formally state:

\noindent
{\bf The optional sampling theorem.}
 Suppose $M(t),\; t \geq 0$ is a martingale with respect to $\eta(t)$ and $\tau$ is a bounded (with probability one) stopping time  with respect to $\eta(t)$, then $\expt[M(\tau)]=\expt[M(0)]$.

Heuristically, the optional sampling theorem states that different realizations of a martingale $M(t)$ stopped at random times average out to constitute a fair game. The crucial aspect is that the stationarity of the expected value holds even for random (stopping) times, {\sfr including our first passage time $\tau$.   As we will see in \S \ref{subsec:app-ddm}, this stationarity enables us to calculate analytic expressions for first passage time properties by finding appropriate martingales.}

A helpful example is to consider the standard Wiener process with initial position at $X(0) = x_0$, and absorbing thresholds at $\pm \zeta${\sfr , with $x_0$ between $\pm \zeta$}. In this case $X(t)$ is itself a martingale, the first passage time $\tau$ is a stopping time, and the optional sampling theorem says that $\expt[ X(0)] = \expt[X(\tau)]$.  The expectation on the left hand side is simply the average of the initial distribution of $X(0)$ which is the number $x_0$. The right hand side is more interesting: $X(\tau)$ is the random value of $X(t)$ at the random decision time $\tau$, the instant $X(t)$ crosses $+\zeta$ or $-\zeta$.  Thus $X(\tau)$ attains one of two possible values, $+\zeta$ or $-\zeta$, and $\expt[X(\tau)]$ heuristically resembles an average over all $+\zeta$'s and $-\zeta$'s corresponding to sample paths starting at $x_0$ and diffusing until they hit either $+\zeta$ or $-\zeta$ at $\tau$.  The optional sampling theorem says that this average of $+\zeta$'s and $-\zeta$'s, over all such sample paths, ends up being equal to the number $x_0$.

\subsection{Applications to the single-stage {\sfr model}}\label{subsec:app-ddm}
The optional sampling theorem is a powerful mathematical tool for decision-making models that associate decisions and decision times with a diffusion processes crossing a threshold. Ratcliff's Diffusion Decision Model \citep{Rat78, Ratcliff2008}, the leaky competing accumulator model \citep{Usher2001}, and the EZ diffusion model \citep{Wagenmakers2007}, are popular examples of such models. The optional sampling theorem reduces the problem of computing analytic expressions for the statistics of the first passage times to identifying appropriate martingales.
In this section we illustrate the flavor of such calculations for the single-stage {\sfr model from} \eqref{eq:pure-ddm} in \S \ref{subsec:ddm}. Recall that the decision time $\tau$ is defined by
$\tau = \inf \{ t>0 | x(t) \notin (-\zeta,\zeta)\}$. Throughout this section, we introduce $\theta = (\mu,\sigma,\zeta)$ to slightly condense the notation when desired.

We first compute  $\prob[x(\tau) =-\zeta]$, the probability of hitting the lower threshold. First, we let $s = \mu/\sigma^2$, the ratio of the drift parameter to the squared diffusion parameter (i.e. signal to noise). Recall from \S\ref{subsec:martingales} that for {\sfr $\mu \neq 0$}, $\exp(-2 s X)$ is a martingale. Applying the optional sampling theorem, we get 
\[
\exp(-2 s x_0) =\expt[ \exp(-2s x(\tau))] = \prob[x(\tau)=\zeta] \expt[ \exp(-2s \zeta)] + \prob[x(\tau)=-\zeta] \expt[ \exp(2s \zeta)].
\]
Substituting $\prob[x(\tau)=\zeta] = 1- \prob[x(\tau)=-\zeta]$ and solving for $\prob[x(\tau)=-\zeta]$, we obtain a closed form expression
\begin{equation*}
  \prob[x(\tau) = -\zeta]  =
 \frac{\exp({-2 s x_0}) - \exp({- 2 s \zeta})}{\exp({2s \zeta}) - \exp({- 2 s \zeta})}.
\end{equation*}
Similarly, for $\mu=0$, we note that $X(t)$ is a martingale. Applying optional sampling theorem, we get
\[
x_0 = \expt[X(\tau )] = \zeta \prob[x(\tau)=\zeta] -\zeta \prob[x(\tau)=-\zeta],
\]
and following the same argument we obtain $  \prob[x(\tau ) = -\zeta]  ={(\zeta-x_0)}/{2\zeta}.$ In summary, we get
\begin{equation}\label{eq:er-ddm} 
 \er (x_0,\theta) := 1- P^+(x_0,\theta) :=  \prob[x(\tau ) = -\zeta] =\begin{cases}\displaystyle
 \frac{\exp({-2 s x_0}) - \exp({- 2 s \zeta})}{\exp({2s \zeta}) - \exp({- 2 s \zeta})}, & \text{if } \mu \ne 0, \\\displaystyle
 \frac{\zeta-x_0}{2\zeta}, & \text{if } \mu = 0,
 \end{cases}
 \end{equation}
where $P^\pm(x_0,\theta)$ is the probability of hitting the upper and the lower threshold, respectively. 

To compute the expected decision time $\expt [ \tau]$ for $\mu \ne 0$, recall from \S\ref{subsec:martingales} that $X(t) - \mu t$ is a martingale. For $\mu \ne 0$, applying the optional sampling theorem yields
\[
x_0 = \expt[(\zeta - \mu \tau)] \prob[x(\tau)= \zeta] + \expt[(-\zeta -\mu\tau)] \prob[x(\tau)= -\zeta]. 
\]
Solving for $\expt[\tau]$, we get $\expt[\tau] = (1-2 \prob[x(\tau)=-\zeta])\zeta/\mu$. When $\mu=0$, recall from \S\ref{subsec:martingales} that $X(t) - \sigma_1^2 t$ is a martingale, and the same argument as above yields $\expt[\tau]= (\zeta^2 -x_0^2)/\sigma_1^2$. In summary, the mean decision time  $\dt(x_0,\theta)$ is given by
\begin{equation}\label{eq:mdt-ddm}
 \dt(x_0,\theta) =
  \expt[\tau] = 
  \begin{cases} \displaystyle
    \frac{(1-2 P^-(x_0, \theta))\zeta -x_0}{\mu}, &\quad \text{if } \mu \ne 0, \\ \displaystyle
    \frac{\zeta^2 -x_0^2}{\sigma^2}, &\quad \text{if } \mu =0.
  \end{cases}
\end{equation}

We also wish to find {\sfr $\tau$'s Laplace transform, or moment generating function\footnote{The moment-generating function {\sfr (technically, the two-sided Laplace transform)} of a random variable $X$ is $\phi_X(\alpha) = \expt [\exp({\alpha X})]$, a function of $\alpha \in \mc A \subset \real $. It is often of interest because it specifies the probability distribution of $X$, and can be used to obtain the moments of $X$.}, $\expt [ \exp ({\alpha \tau})]$.}
We remember from \S\ref{subsec:martingales} that $\exp(\lambda X(t) -\lambda \mu t - \lambda^2 \sigma^2 t /2)$ is a martingale, and choose $\lambda$ so that the coefficient of $t$ becomes $-\alpha$, i.e., $\lambda$ solves the equation $\frac{1}{2}\sigma^2 \lambda^2 + \mu \lambda - \alpha=0$. The two solutions to this equation are 
\[
\lambda_1 = \frac{-\mu -\sqrt{\mu + 2\alpha \sigma^2}}{\sigma^2}, \quad \text{and}  \quad 
\lambda_2 = \frac{-\mu +\sqrt{\mu + 2\alpha \sigma^2}}{\sigma^2}.
\]
Applying the optional sampling theorem, we obtain a pair of equations,
\begin{align*}
\exp(\lambda_1 x_0) &= \exp({\lambda_1 \zeta})\expt[\exp(-\alpha \tau ) 1(x(\tau)=\zeta)] + \exp({-\lambda_1 \zeta})\expt[\exp(-\alpha \tau ) 1(x(\tau)=-\zeta)] \\ 
\exp(\lambda_2 x_0) &= \exp({\lambda_2 \zeta})\expt[\exp(-\alpha \tau ) 1(x(\tau)=\zeta)] + \exp({-\lambda_2 \zeta})\expt[\exp(-\alpha \tau ) 1(x(\tau)=-\zeta)],
\end{align*}
which we can solve simultaneously for $\expt[\exp(-\alpha \tau ) 1(x(\tau)= \pm \zeta)]$ to obtain
\[
\expt[\exp({-\alpha \tau }) \bs 1(x(\tau)= \pm \zeta)] = \exp\Big( \frac{\mu (\pm \zeta-x_0)}{\sigma^2}\Big) \frac{ \sinh (\frac{(\zeta \pm x_0) \sqrt{2\alpha \sigma^2 + \mu^2} }{\sigma^2}) }{\sinh (\frac{2\zeta\sqrt{2\alpha \sigma^2 + \mu^2}}{\sigma^2})}.
\]
Thus, the moment generating function of the decision time is
\begin{equation}\label{eq:laplace-ddm}
\expt[\exp({-\alpha \tau })] = \exp\Big( \frac{\mu ( \zeta-x_0)}{\sigma^2}\Big) \frac{ \sinh (\frac{(\zeta + x_0) \sqrt{2\alpha \sigma^2 + \mu^2} }{\sigma^2}) }{\sinh (\frac{2\zeta\sqrt{2\alpha \sigma^2 + \mu^2}}{\sigma^2})} + \exp\Big( -\frac{\mu ( \zeta +x_0)}{\sigma^2}\Big) \frac{ \sinh (\frac{(\zeta - x_0) \sqrt{2\alpha \sigma^2 + \mu^2} }{\sigma^2}) }{\sinh (\frac{2\zeta\sqrt{2\alpha \sigma^2 + \mu^2}}{\sigma^2})}.
\end{equation}
As a byproduct, we also get the Laplace transform of conditional decision times: 
\begin{align}\label{eq:laplace-ddm-correct}
\expt[\exp({-\alpha \tau}) | x(\tau)=\zeta] &= \frac{\exp( \frac{\mu (\zeta-x_0)}{\sigma^2})}{P^+(x_0, \theta)} \frac{ \sinh (\frac{(\zeta+x_0) \sqrt{2\alpha \sigma^2 + \mu^2} }{\sigma^2}) }{\sinh (\frac{2\zeta\sqrt{2\alpha \sigma^2 + \mu^2}}{\sigma^2})}\\\label{eq:laplace-ddm-incorrect}
 \expt[\exp({-\alpha \tau}) |x(\tau)=-\zeta] & = \frac{\exp(- \frac{\mu (\zeta + x_0)}{\sigma^2})}{P^-(x_0, \theta)} \frac{ \sinh (\frac{(\zeta-x_0) \sqrt{2\alpha \sigma^2 + \mu^2} }{\sigma^2}) }{\sinh (\frac{2\zeta\sqrt{2\alpha \sigma^2 + \mu^2}}{\sigma^2})}.
\end{align}

The derivatives of the Laplace transform yield moments of decision time (see ~\cite{VS-PH-PS:14-arxiv} for detailed derivation of conditional and unconditional moments of decision time using Laplace transforms). Here, we focus on expressions for conditional expected (mean) decision times that are the derivative of the Laplace transform with respect to $-\alpha$ computed at $\alpha=0$.
The expected decision time conditioned on hitting the upper and lower boundaries are denoted by $\dt^+$ and $\dt^-$, and may be computed by differentiating \eqref{eq:laplace-ddm-correct} and \eqref{eq:laplace-ddm-incorrect}:
\begin{align} \label{eq:dtplus-ddm}
  \dt^+(x_0, \theta) &= \frac{\widehat{\dt}^+(x_0,\theta)}{P^+(x_0, \theta)} =
  \begin{cases}
 \frac{2\zeta}{\mu} \coth( {2 s \zeta} ) - \frac{\zeta+x_0}{\mu} \coth ( { s (\zeta +x_0)}), & \text{if } \mu \ne 0 ,\\
 \big( \frac{4\zeta^2}{3 \sigma^2} - \frac{(\zeta+x_0)^2}{3 \sigma^2}\big), & \text{if } \mu =0;
 \end{cases} \\ \label{eq:dtminus-ddm}
  \dt^-(x_0, \theta) &= \frac{\widehat{\dt}^-(x_0,\theta)}{P^-(x_0, \theta)} =
  \begin{cases}
 \frac{2\zeta}{\mu} \coth( {2 s \zeta} ) - \frac{\zeta-x_0}{\mu} \coth ( { s (\zeta -x_0)}), & \text{if } \mu \ne 0 ,\\
\big( \frac{4\zeta^2}{3 \sigma^2} - \frac{(\zeta-x_0)^2}{3 \sigma^2}\big), & \text{if } \mu =0,
 \end{cases}
\end{align}
where $\widehat{\dt}^\pm = \expt[\tau \bs 1(x(\tau) = \pm \zeta)]$ and $\bs 1(\cdot)$ is the indicator function. 
We again note, just as with $\prob [x(\tau)=-\zeta]$,  $\dt^\pm$ also depend on the underlying parameters $\mu, x_0, \zeta$, and $\sigma$.

We now compute $\tau$'s first passage time density $f(t;x_0,\mu,\sigma,\zeta)$, i.e., the probability density function of the decision time.  {\sfr This} amounts to calculating the inverse Laplace transform of ~\eqref{eq:laplace-ddm}. In this case, the inverse Laplace transform needs to be expressed as an infinite series (see~\cite{Lin1998} for a detailed derivation):
\begin{multline}\label{eq:fpt-density-ddm}
f(t;x_0,\theta) = \frac{d}{dt}\prob [ \tau \le t ] =  \exp\Big({-\frac{a^2 t}{2\sigma^2}}\Big) \bigg(\exp({-s(\zeta + x_0)}) \vartheta\Big(t; \frac{\zeta-x_0}{\sigma}, \frac{2\zeta}{\sigma}\Big) \\
+ \exp({s(\zeta-x_0)})  \vartheta\Big(t; \frac{\zeta+ x_0}{\sigma}, \frac{2\zeta}{\sigma}\Big) \bigg), 
\end{multline}
where $\theta=(a,\sigma,\zeta)$ , and $\vartheta(t; u,v)$ is a function~\cite[pp. 451]{Borodin2002} defined by
\[
\vartheta(t; u,v) = \sum_{k=-\infty}^{+\infty} \frac{v-u + 2 kv}{\sqrt{2 \pi} t^{3/2}} \exp\Big({-\frac{(v-u+2kv)^2}{2t}}\Big), \quad u<v. 
\]

Similarly, the first passage time density conditioned on a particular decision is given by
\begin{align} \label{eq:cond-fpt-ddm-1}
\frac{d}{dt}\prob [\tau \le t | x(\tau) =\zeta] &=\frac{f^{+}(t;x_0,\theta)}{P^+(x_0, \theta)}   = \frac{\exp({-\frac{\mu^2 t}{2\sigma^2}+ s(\zeta-x_0)})}{P^+(x_0, \theta)}  \vartheta\Big(t; \frac{\zeta+ x_0}{\sigma}, \frac{2\zeta}{\sigma}\Big) \\ \label{eq:cond-fpt-ddm-2}
\frac{d}{dt}\prob[\tau \le  t | x(\tau)=-\zeta] & =\frac{f^{-}(t;x_0,\theta)}{P^-(x_0, \theta)} = \frac{\exp({-\frac{\mu^2 t}{2\sigma^2}-s(\zeta + x_0)})}{P^-(x_0, \theta)}  \vartheta\Big(t; \frac{\zeta-x_0}{\sigma}, \frac{2\zeta}{\sigma}\Big),
\end{align}
where $f^{\pm}(t;x_0,\theta) = \frac{d}{dt} \prob [\tau \leq t\; \& \;x(\tau) = \pm \zeta]$, i.e., $f^{\pm}(t;x_0,\theta) dt$ is the probability of the event $\tau \in [t, t+ dt)$ and $x(\tau) =\pm \zeta$.  Note that $f$ defined in~\eqref{eq:fpt-density-ddm} is the sum of $f^+$ and $f^-$.

{Alternate derivations}
for the hitting probabilities, mean decision times, and FPT densities may be found in the decision making literature \citep{Ratcliff2004,Bogacz2006,Navarro2009}.  
It is worth noting that the infinite series solution for the FPT density given in~\eqref{eq:fpt-density-ddm} is equivalent to the small-time representations for the FPT analyzed in~(\cite{Navarro2009}~and~\cite{Blurton2012}). For completeness, we provide the alternative expression for density in~\ref{app:fpt-density}.

\section{Analysis of the two-stage {\sfr model}}
\label{sec:2-ddm}
In this section, we use the tools developed in \S \ref{sec:review} in order to analyze the two-stage {\sfr process} originally presented and analyzed in \cite{Ratcliff1980}.
While our calculations lead to equivalent formulas for the first passage time densities, a martingale argument provides us with additional closed form expressions for {\sfr the} probability of hitting a particular threshold and expected decision times.
Computations of these FPT statistics using {\sfr the results of \cite{Ratcliff1980} requires} numerical integration of the FPT density, which our formulas now avoid.

We may explicitly write the two-stage {\sfr model} as
\begin{equation} \label{eq:2-ddm}
dx(t) = \mu(t) \, dt + \sigma (t)dW(t),\,\,\, x(t_0) = x_0,
\end{equation}
where 
\begin{align*}
\mu(t) = \begin{cases}
\mu_1, &  \text{ for }  0 \leq t < t_1, \\
\mu_2, &  \text{ for }   t \ge  t_1,
\end{cases}  \quad \text{and} \quad \sigma(t) = \begin{cases}
\sigma_1, &  \text{ for }  0 \leq t < t_1, \\
\sigma_2, &  \text{ for }   t \ge  t_1.
\end{cases} 
\end{align*}
As before, the decision time $\tau = \inf \{ t>0 | x(t) \notin (-\zeta,\zeta)\}$ is the first passage time with respect to boundaries at $\pm \zeta$. Let $\theta_i =(\mu_i, \sigma_i, t_i, \zeta)$ and $s_i =\mu_i/\sigma_i^2$, for $i \in \{1,2\}$.

We interpret the first stage as a single-stage {\sfr model} with a deadline at $t_1$. For a single-stage {\sfr model} with thresholds $\pm \zeta$ and a deadline $t_1$, the joint density {\sfr $\supscr{g}{ddln}(x,\tau; x_0,t_1, \theta_1)$} of the evidence $x(t_1)$ and the event $\tau \ge  t_1$ is given by {\sfr \citep{Douady1999,Durrett2010}:}
{\sfr
\begin{align}
  \supscr{g}{ddln}(x,\tau ; x_0, t_1, \theta_1) &= \frac{d}{d x} \prob[x(t_1) \le x \; \& \; \tau \ge t_1] \nonumber\\
  &= \frac{\boldsymbol 1[x \in (-\zeta,\zeta) ]}{\sqrt{2 \pi t \sigma_1^2}}
  \exp \left(\frac{-\mu_1^2 t + 2 \mu_1 (x-x_0)}{2\sigma_1^2} \right) \times \nonumber\\  \label{eq:dist-deadline-fixed}
  &  \sum_{n=-\infty}^{\infty} \left[ \exp \left(\frac{-(x-x_0+ 4n\zeta)^2}{2\sigma_1^2 t}\right) -  \exp \left(\frac{-(2\zeta-x -x_0 + 4n\zeta)^2}{2\sigma_1^2 t}\right) \right]. 
\end{align}
}
Here superscript ``ddln" refers to the {\sfr deadline.  
  $\supscr{g}{ddln}$} may then be used to determine the FPT distribution by integrating it over the range of $x$. More importantly, dividing $\supscr{g}{ddln}$ by $\prob [ \tau \ge t_1 ]$ yields the conditional density on the evidence $x(t_1)$ conditioned on no decision until time $t_1$.

\subsection{Probability of hitting the lower threshold}
In trying to compute $\prob [x(\tau)=-\zeta]$, we view the two-stage process as {\sfr two single-stage processes in sequence}. Let $\tau_1$ be the first passage time for the first stage by itself (without any deadline at $t_1$) and define the random time $  \hat \tau_1 = \min\{\tau_1, t_1\}$ which is a stopping time. Applying the optional sampling theorem to the exponential martingale for the first stage gives us
\begin{align*}
&\quad \exp({-2 s_1 x_0}) = \expt[\exp({-2s_1 x(\hat \tau_1)})] \\
&= \expt[\exp({-2s_1 x(\tau_1)}) | \tau_1 \le t_1] \prob [ \tau_1 \le t_1 ] + 
\expt[\exp({-2s_1 x(t_1)})| \tau_1 > t_1] \prob [ \tau_1 > t_1 ] \\
& = \left(\exp({-2s_1 \zeta}) \prob [x(\tau_1)=\zeta|\tau_1\le t_1]+ \exp({2s_1 \zeta}) \prob[x(\tau_1)=-\zeta|\tau_1\le t_1] \right) \prob [ \tau_1 \le t_1 ] \\
& \qquad \qquad  \qquad \qquad  \qquad \qquad \qquad \qquad  \qquad \qquad +  \expt[\exp({-2s_1 x(t_1)}) \bs 1(\tau_1 > t_1)] \\
& = \left(\exp({-2s_1 \zeta}) + ( \exp({2s_1 \zeta}) - \exp({-2s_1 \zeta})) \prob[x(\tau_1)=-\zeta|\tau_1\le t_1]\right)\prob [ \tau_1 \le t_1 ] \\
&\qquad \qquad  \qquad \qquad  \qquad \qquad \qquad \qquad  \qquad \qquad
+  \expt[\exp({-2s_1 x(t_1)}) \bs 1(\tau_1 > t_1)],
\end{align*} 
which yields
\begin{align*}
\prob[x(\tau_1)=-\zeta \; \& \; \tau_1\le t_1] = 
\frac{\exp({-2 s_1 x_0}) -  \expt[\exp({-2s_1 x(t_1)}) \bs 1(\tau_1 > t_1)] - \exp({-2s_1 \zeta}) \prob [ \tau_1 \le t_1 ]}{ \exp({2s_1 \zeta}) - \exp({-2s_1 \zeta})}.
\end{align*}
In this expression, both $\expt[\exp({-2s_1 x(t_1)}) \bs 1(\tau_1 > t_1)]$ and $ \prob [ \tau_1 \le t_1 ]$ may be obtained from~\eqref{eq:dist-deadline-fixed}.

The {\sfr second stage} is another single-stage {\sfr process}, this time starting at time $t_1$ with a random initial condition $x(t_1)$ with distribution given by~\eqref{eq:dist-deadline-fixed}. 
Computing the expected value of the standard lower threshold hitting probability~\eqref{eq:er-ddm} with respect to the random initial condition $X_1$, i.e., $x(t_1)$ conditioned on $\tau_1 > t_1$ we obtain
\[
\prob [x(\tau)=-\zeta | \tau_1 > t_1] = \frac{\expt[\exp(-2s_2 x(t_1))|\tau_1 > t_1] -\expt(-2s_2 \zeta)}{\exp(2s_2 \zeta) - \exp(-2s_2 \zeta)}.
\]

The $\expt[\exp(-2s_2 x(t_1))|\tau_1 > t_1]$ term can be readily computed from~\eqref{eq:dist-deadline-fixed}. Combining the previous two conditional expressions allows us to obtain
\[
\prob[x(\tau)=-\zeta] = \prob[x(\tau_1)=-\zeta \; \& \; \tau_1\le t_1]  + \prob[x(\tau)=-\zeta | \tau_1 > t_1]  \prob [ \tau_1>t_1 ]. 
\]
This probability depends on all of the parameters of the two-stage model: $x_0,\mu_1,\sigma_1, t_1, \mu_2,\sigma_2,\zeta$.

\subsection{Expected decision time}
To compute the expected decision time, we apply the optional sampling theorem to the martingale $x(t) - \mu_1 t$ with the stopping time $\hat \tau_1$ from the previous section to obtain

\begin{align*}
x_0 &= \expt[x(\hat \tau_1) -\mu_1 \hat \tau_1] = 
 \expt[x(\tau_1) -\mu_1  \tau_1| \tau_1 \le t_1] \prob [ \tau_1 \le t_1 ] +  \expt[x(t_1) -\mu_1 t_1| \tau_1 > t_1] \prob [ \tau_1 > t_1 ]\\
 & = (\zeta (1 -2 \prob[x(\tau_1)=-\zeta |\tau_1 \le t_1 ])- \mu_1 \expt[\tau_1|\tau_1 \le t_1]) \prob [ \tau_1 \le t_1 ] +( \expt[x(t_1)|\tau_1 >t_1] - \mu_1 t_1) \prob [ \tau_1 >t_1].
\end{align*}
Solving for $\expt[\tau \bs 1(\tau_1 \le t_1)]$ we obtain
\[
\expt[\tau\bs 1(\tau_1 \le t_1)] = \frac{(\prob [ \tau_1\le t_1 ] -2 \prob[x(\tau_1)=-\zeta \; \& \; \tau_1\le t_1 ] ) \zeta - x_0 +(\expt[x(t_1)\bs 1(\tau_1 >t_1)] - \mu_1 t_1 \prob [ \tau_1 >t_1 ]) }{\mu_1}.
\]
Much like in the previous section, we observe that the second stage is similar to a single-stage {\sfr process} starting at time $t_1$ with a random initial condition determined by~\eqref{eq:dist-deadline-fixed}. Thus, the associated expected first passage time is
\[
\expt[\tau | \tau_1 >t_1] = t_1 + \frac{(1- 2\prob[x(\tau)=-\zeta | \tau > t_1])\zeta - \expt[x(t_1)|\tau_1 >t_1]}{\mu_2}.
\]
Combining the above expressions, we obtain
\begin{multline*}
\expt[\tau] = \frac{(\prob [ \tau_1\le t_1 ] -2 \prob[x(\tau_1)=-\zeta \; \& \; \tau_1\le t_1 ] )\zeta - x_0 +(\expt[x(t_1)\bs 1(\tau_1 >t_1)] }{\mu_1}  \\
+ \frac{(1- 2\prob [x(\tau)=-\zeta | \tau > t_1]) \prob [ \tau_1 > t_1 ]\zeta - \expt[x(t_1) \bs 1(\tau_1 >t_1)]}{\mu_2}.
\end{multline*}

\subsection{First passage time density}
We now compute the first passage time probability density function. 
Let $F(t; \theta, x_0)$ be the cumulative distribution function of the decision time for the single-stage process~\eqref{eq:pure-ddm} obtained by integrating~\eqref{eq:fpt-density-ddm}.
For $t \le t_1$, the {\sfr two-stage model} is identical to the {\sfr first-stage model} and the first passage time distribution is $F(t; \theta_1, x_0)$. For $t>t_1$, the first passage time distribution is 
\[
F(t; \theta_1, \theta_2, x_0)= F(t_1; \theta_1, x_0) + \expt[F(t- t_1; \theta_2, x(t_1)) \bs 1(\tau_1 \ge t_1)],
\]
i.e., the distribution function corresponds to trajectories that {\sfr reach} threshold before {\sfr $t_1$} and trajectories that {\sfr reach} threshold between $t_1$ and $t$. The latter trajectories can be modeled as trajectories of a {\sfr single stage process} starting at time $t_1$ with stochastic initial condition. The stochastic initial condition leads to the expectation operator on the second term. The distribution of decision time conditioned on particular decisions can be computed analogously. 

\section{Analysis of the {\sfr multistage model}}
\label{sec:performance-metrics}
In this section we derive first passage time (FPT) properties of the
{\sfr multistage process} defined in \S\ref{sec:msdm} using an
approach similar to that employed throughout \S\ref{sec:2-ddm}. The
{\sfr model} is viewed as {\sfr $n$ modified processes in sequence}
in which for each {\sfr stage,} the initial condition is a random
variable and only the decisions made before a deadline are
considered. For the $i$-th stage {\sfr process} with a known
distribution of initial condition $X_{i-1}$, we derive properties of
the FPT conditioned on a decision before the deadline $t_i$, along
with the distribution of $X_i$ for $i\in \until{n}$. This latter
distribution, more precisely, is the distribution of $x(t_i)$
conditioned on the FPT for the $i$-th {\sfr stage being} greater than
$t_i$. We then use these properties sequentially for $i \in \until{n}$
to determine the FPT properties during each stage. Finally, we
aggregate FPT properties at each stage to compute FPT properties for
the {\sfr whole multistage process}.  Our calculations have features
similar to an idea of \cite{Diederich2006}, who first proposed that
the bias (i.e. initial condition) of a stage may have a time
dimension. In some sense our formulas below elucidate how previous
temporal stages of processing affect the bias of future stages.

The extension of the FPT distribution computation {\sfr from the two-stage model}~\citep{Ratcliff1980} to {\sfr the multistage model} requires careful computation of expressions similar to~\eqref{eq:dist-deadline-fixed} at the end of each stage. Also, in contrast to~\cite{Ratcliff1980}, as in \S\ref{sec:2-ddm}, our martingale based approach allows direct computation of probability of hitting a particular threshold and expected decision times. {As in the previous section, this avoids} integration of {\sfr the} first passage time density to compute these quantities{\sfr. }

Throughout this section we use the following notations:
\begin{itemize}
\item $\tau = \inf \{ t>0 | x(t) \notin (-\zeta,\zeta)\}$, the first passage time through either threshold for the entire multistage process;
\item $\tau_i = \tau |  \tau > t_{i-1}$, the first passage time for the $i$-th {\sfr stage }~\eqref{eq:ddm-i} without any deadline;
\item $\theta_i = (\mu_i,\sigma_i, t_i, \zeta)$, and $\theta_{1:\ell} = (\mu_1,\ldots,\mu_\ell,\sigma_1,\ldots,\sigma_\ell,t_1, \ldots, t_\ell, \zeta)$ representing the parameters for the $i$-th stage and stages $1,\ldots, \ell$, respectively;
\item $s_i = \mu_i/\sigma_i^2$, the $i$-th stage ratio of signal to squared noise.
\end{itemize}

Here we are concerned with computation of FPT properties and allow $\theta_{1:n}$ to be free parameters. For scenarios such as estimation of parameters, to ensure identifiability of the parameters, all diffusion rates may be  set equal to unity. However, such cases are beyond the scope of this manuscript and we do not discuss these issues here. 

\subsection{FPT properties of the $i$-th stage}\label{subsec:i-ddm}
For the {\sfr $i$-th stage,} the initial condition $X_{i-1}$ is a random variable and only decisions made before the deadline $t_i$ are relevant. The analysis thus focuses on the random variable $X_{i-1}$ and the random time $\tau_i$. Conditioned on a realization of $X_{i-1}$, the density of $X_i$ can be computed using~\eqref{eq:dist-deadline-fixed}. If the density of $X_{i-1}$ is known, then the unconditional density of $X_i$ can be obtained by computing the expected value of the conditional density of $X_i$ with respect to $X_{i-1}$. Since the density of $X_0$ is known, this procedure can be recursively applied to obtain densities of $X_{i-1}$, for each $i\in \until{n}$.  Formally, the joint density {\sfr  $\supscr{g}{ddln}_i(x, \tau ; x_0, \theta_{1:i})$} of  the evidence $x(t_i)$ and the event $\tau \geq t_i$ is
{\sfr
\begin{align}
\supscr{g}{ddln}_i(x, \tau_i ;x_0, t_i , \theta_{1:i}) &=  \frac{d}{d x} \prob[ x(t_i) \le x  \; \& \; \tau_i > t_i]
 \label{eq:dist-deadline}
= \expt_{X_{i-1}}[\supscr{g}{ddln}(x, \tau_i -t_{i-1}; X_{i-1},t_i-t_{i-1}, \theta_i)], 
\end{align}
}
where $\expt_{X_{i-1}}[\cdot]$ denotes the expected value with respect to $X_{i-1}$. Note that 
\begin{align*}
 \prob[ x(t_i) \le x  \; \& \; \tau_i > t_i] = \prob[ x(t_i) \le x  \; | \; \tau_i > t_i] \prob [ \tau_i > t_i  ]=  \prob [  X_i \le x  ] \prob [ \tau_i > t_i  ].
\end{align*}
Thus, the density of $X_i$, i.e., $x(t_i)$ conditioned on $\tau_i > t_i$  is determined by dividing $\supscr{g}{ddln}_i$ by $\prob [ \tau_i > t_i ]$ which can be computed by integrating $\supscr{g}{ddln}_i$ over the range of $x(t_i)$. Note that the parameters in $\supscr{g}{ddln}_i$ are $x_0$ and $\theta_{1:i}$; this highlights the fact that the distribution of $X_i$ depends on all previous stages. 

Similarly, the FPT density for the $i$-th {\sfr stage} conditioned on a realization of $X_{i-1}$ can be computed using~\eqref{eq:fpt-density-ddm}, and the unconditional density can be obtained by computing the expected value of the conditional density with respect to $X_{i-1}$:
\begin{equation} \label{eq:fpt-density}
f_i(t;x_0,\theta_{1:i}) = \frac{d}{dt}\prob [  \tau_i  \le  t ] = \expt_{X_{i-1}} [f(t-t_{i-1};X_{i-1},\theta_i) ],  
\end{equation}
where $t > t_{i-1}$. The cumulative distribution function $F_i(t;x_0,\theta_{1:i}) = \prob [ \tau_i \leq t ]$   is obtained by integrating $f_i(t;x_0,\theta_{1:i})$. Note that every trajectory crossing the decision threshold before $t_i$ does so irrespective of the deadline at $t_i$. Thus, the expression for density $f_i$ does not depend on $t_i$. 

To conclude this section, we state equations for computing hitting times, mean decision times, and first passage time densities, conditional on a response during the $i$-th stage. The derivation of these expressions are found in~\ref{app:ddm-i}.
\begin{enumerate}[(i)]
\item The probability of hitting the lower threshold given that a response is made during the $i$-th stage, denoted by $\er_i(x_0, \theta_{1:i}) := 1- P^+_i(x_0, \theta_{1:i}) :=\prob [ x(\tau) =  -\zeta| t_{i-1}<  \tau \le  t_i] $, is given by
\begin{align}\label{eq:error-rate-deadline}
\resizebox{0.85\textwidth}{!}{$
\er_i= \begin{cases} \displaystyle
\frac{\expt_{X_{i-1}} [\exp(-2s_i X_{i-1})] - \expt_{X_i} [\exp(-2s_i X_i)] \prob [ \tau_i > t_i ]  - \exp({-2s_i \zeta}) \prob [ \tau_i \le t_i ]}{( \exp({ 2s_i \zeta}) - \exp({-2s_i \zeta}) ) \prob [ \tau_i \le  t_i ]}, & \text{ if } \mu_i \ne 0, \\ 
&\\
 \displaystyle
\frac{1}{2} - \frac{(\expt_{X_{i-1}}[X_{i-1}] - \expt_{X_i}[X_i ] \prob [ \tau_i >t_i ])}{2 \zeta \prob [ \tau_i \le t_i ]}, & \text{ if } \mu_i=0.
\end{cases}$}
\end{align}
$P_i^+(x_0, \theta_{1:i})$ {\sfr is the probability} of hitting the upper threshold during the $i$-th stage. These expression depends on all of the {\sfr model} parameters up to and including the $i$-th stage (i.e. $P_i^- := P_i^- (x_0,\theta_{1:i})$ and $P_i^+ := P_i^+(x_0,\theta_{1:i})$ both depend on $x_0$ and $\theta_{1:i}$ ) and will be used in subsequent calculations.

\item The joint FPT density for the $i$-th stage {\sfr process} and a given upper/lower response, denoted by $f_i^\pm(t;x_0,\theta_{1:i}) = \frac{d}{dt} \prob [\tau_i \leq t \; \& \; x(\tau_i) = \pm \zeta]$, is given by
\begin{align}
f^+_i(t;x_0,\theta_{1:i}) &= \frac{d}{dt}\prob [ \tau_i  \le  t \; \& \: x(\tau_i) = \zeta ]  = \expt_{X_{i-1}} [f^+(t-t_{i-1};X_{i-1},\theta_i) ] \label{eq:fpt-density-correct} \\
f^-_i(t;x_0,\theta_{1:i}) &=\frac{d}{dt}\prob [ \tau_i  \le  t \; \& \: x(\tau_i) = -\zeta] = \expt_{X_{i-1}} [ f^-(t-t_{i-1};X_{i-1},\theta_i) ], \label{eq:fpt-density-incorrect}
\end{align}
where the functions $f^\pm(t;x_0,\theta_i)$ are taken from~\eqref{eq:cond-fpt-ddm-1}~and~\eqref{eq:cond-fpt-ddm-2}.
Again, these expressions depend on all of the {\sfr multistage model} parameters up to and including the $i$-th stage (i.e. both $f^\pm_i(t;x_0,\theta_{1:i}) $ depend on $x_0$ and $\theta_{1:i}$). 

\item The mean decision time given a response during stage $i$, denoted by $\dt _i(x_0,\theta_{1:i})$, is given by
\begin{align} 
&\dt _i(x_0,\theta_{1:i}) = \expt[\tau_i | \tau_i \le t_i] \nonumber\\ \label{eq:decision-time-deadline}
&= 
\resizebox{0.85\textwidth}{!}{$
\begin{cases} \displaystyle 
t_{i-1} +   \frac{(1-2 \er_i)\zeta  \prob [ \tau_i \le t_i ]  - \expt_{X_{i-1}} [X_{i-1}] + \expt_{X_i}[X_i] \prob [ \tau_i  > t_i ] - \mu_i (t_i -t_{i-1}) \prob [ \tau_i > t_i ]}{\mu_i \prob [ \tau_i \le t_i ]}, & \text{ if } \mu_i \ne 0, \\  
&\\
\displaystyle
t_{i-1} +  \frac{\zeta^2 \prob [ \tau_i \le t_{i} ] - \expt_{X_{i-1}}[X_{i-1}^2] + \expt_{X_i}[X_i^2] \prob [ \tau_i > t_i ] - \sigma_i^2 (t_i - t_{i-1})\prob [ \tau_i > t_i ]}{\sigma_i^2 \prob [ \tau_i \le t_i ]}, & \text{ if } \mu_i=0.
 \end{cases}$}
\end{align}

\item The mean decision time conditioned on a given upper/lower response made during the $i$-th stage, denoted by $\dt _i^\pm(x_0, \theta_{1:i})$, is given by
  \begin{align}
    \dt_i^+(x_0, \theta_{1:i}) =& \expt[\tau_i | x(\tau_i) =\zeta  \; \&\; \tau_i \le t_i] \nonumber \\
    =& \frac{\widehat{\dt}_i^+(x_0,\theta_{1:i})}{P^+(x_0, \theta_{1:i}) \prob [ \tau_i \le t_i ]}\nonumber \\
    =&\frac{1}{{P^+(x_0, \theta_{1:i}) \prob [ \tau_i \le t_i ]}} 
    \bigg( t_{i-1}  \prob [x(\tau_i)=\zeta] + 
    \expt_{X_{i-1}} \left[\widehat{\dt}^+(X_{i-1},\theta_i) \right]  \nonumber\\ \label{eq:decision-time-correct-deadline}
    & \qquad \qquad-  \Big( \expt_{X_i} \left[\widehat{\dt}^+(X_{i},\theta_i) \right]  
    - t_i \prob [x(\tau_i) =\zeta | \tau_i > t_i] \Big)\prob [ \tau_i > t_i ] \bigg)
  \end{align}
  \begin{align}
    \dt_i^-(x_0, \theta_{1:i}) =& \expt[\tau_i | x(\tau_i) =-\zeta  \; \&\; \tau_i \le t_i] \nonumber \\
    =& \frac{\widehat{\dt}^-_i(x_0,\theta_{1:i}) }{\er_i(x_0, \theta_{1:i}) \prob [ \tau_i \leq t_i ]}\nonumber \\ \label{eq:decision-time-incorrect-deadline}
    =& \frac{1}{{\er_i(x_0, \theta_{1:i}) \prob [ \tau_i \le t_i ]}} \bigg( t_{i-1}  \prob[x(\tau_i)=-\zeta] +
\expt_{X_{i-1}} \left[\widehat{\dt}^-(X_{i-1},\theta_i) \right]  \nonumber\\
&\qquad \qquad - \Big(\expt_{X_i} \left[ \widehat{\dt}^-(X_i,\theta_i) \right]  - t_i \prob[x(\tau_i) = -\zeta | \tau_i > t_i]\Big)\prob [ \tau_i > t_i ]\bigg),
  \end{align}
  where
  $\displaystyle \widehat{\dt}_i^\pm (x_0,\theta_{1:i}) = \expt[\tau_i \bs 1(x(\tau_i) = \pm \zeta \; \& \; \tau_i \leq t_i)] $ and
 is calculated using~\eqref{eq:dtplus-ddm}~and~\eqref{eq:dtminus-ddm}, $\prob [x(\tau_i)= \pm \zeta]= \expt_{X_{i-1}}[P^\pm(X_{i-1}, \theta_i)] $, and  $\prob [x(\tau_i) = \pm \zeta \;|\; \tau_i > t_i ] = \expt_{X_i}[P^\pm(X_i, \theta_i)]$  is calculated using~\eqref{eq:er-ddm} and~\eqref{eq:dist-deadline}.
\end{enumerate}

\subsection{FPT properties of the {\sfr multistage model}}
\label{subsec:multistage-ddm}
For a given {\sfr multistage process} \eqref{eq:multistage-ddm} with initial condition $x_0$, we sequentially compute all of the distributions of the initial conditions $X_i$ for each $i \in \until{n}$ using~\eqref{eq:dist-deadline}. Then, we compute the properties of the FPT associated with the $i$-th stage. Finally, the total probability formula aggregates these results into FPT properties of the entire {\sfr multistage model}.  The calculations are contained in~\ref{app:mult-ddm} and are given in terms of earlier formulas. In the following, we omit the arguments of functions whenever it is clear from the context. 
\begin{enumerate}[(i)]
\item Let $t \geq 0$ be given such that $t \in (t_{k-1}, t_k]$ for some $k \in \until{n}$. The FPT distribution for the {\sfr multistage process} is
 \begin{align}
\prob [ \tau \le t ] & = 1 - \prod_{i=1}^{k-1} \prob [ \tau_i > t_i ] +  \prob [ \tau_k \le t ] \prod_{i=1}^{k-1} \prob [ \tau_i > t_i ]. \label{eq:multistage-fpt}
\end{align}
Note that $\prod_{i=1}^{k-1} \prob [ \tau_i > t_i ]= \prob [ \tau > t_{k-1} ]$ and $\prob [ \tau_k \le t ] = \prob [ \tau \le t | \tau > t_{k-1} ]$. 

\item The mean decision time, denoted by $\subscr{\dt}{ms}(x_0, \theta_{1:n})$,  {\sfr for the multistage process is}
\begin{align} 
\subscr{\dt}{ms}(x_0, \theta_{1:n}) =\expt[\tau]  &=  \sum_{i=1}^n \Big(\expt[\tau_i | \tau_i \le  t_i ]  \prob [ \tau_i \le  t_i  ] \prod_{j=1}^{i-1} \prob [ \tau_j > t_i ] \Big) \nonumber\\ \label{eq:multistage-DT}
&=  \sum_{i=1}^n \Big(\dt_i(x_0,\theta_{1:i})  \prob [ \tau_i \le  t_i  ] \prod_{j=1}^{i-1} \prob [ \tau_j > t_i ] \Big).
\end{align}
Put simply, the expected decision time is the sum of the expected decision times for the individual stages ($\dt_i$) weighted by the probability of the decision in each stage ($ \prob [ \tau_i \le  t_i  ] \prod_{j=1}^{i-1} \prob [ \tau_j > t_i ] $).

\item The probability of hitting the lower threshold , denoted by $\subscr{P}{ms}^-(x_0, \theta_{1:n})$, {\sfr is}
\begin{align} \label{eq:multistage-ER}
\subscr{P}{ms}^-(x_0, \theta_{1:n})= 1- \subscr{P}{ms}^+ (x_0, \theta_{1:n}) =\sum_{i=1}^n \Big( \er_i \prob [ \tau_i < t_i  ] \prod_{j=1}^{i-1} \prob [ \tau_j > t_j ] \Big).
\end{align}
This expression is similar to \eqref{eq:multistage-DT}, with the
probability of hitting the lower threshold being the sum of the hitting
probabilities for each stage ($\er_i $) weighted by the probability of
the decision in each stage ($ \prob [ \tau_i \le t_i ]
\prod_{j=1}^{i-1} \prob [ \tau_j > t_i ] $).

\item The mean decision time conditioned on hitting the 
upper/lower threshold is
\begin{align} \label{eq:multistage-DT-correct}
\expt[\tau | x(\tau) =\zeta ]  &=  \frac{1}{\subscr{P}{ms}^+}\sum_{i=1}^n \Big( \expt[\tau_i | x(\tau_i)=\zeta \; \&\; \tau_i \le  t_i ] \;  P_i^+ \; \prob [ \tau_i \le  t_i  ] \prod_{j=1}^{i-1} \prob [ \tau_j > t_i ] \Big) \\ \label{eq:multistage-DT-incorrect}
\expt[\tau | x(\tau) =-\zeta ]  &= \frac{1}{\subscr{P}{ms}^-} \sum_{i=1}^n \Big( \expt[\tau_i | x(\tau_i)=\zeta \; \&\; \tau_i \le  t_i ] \;  \er_i \; \prob [ \tau_i \le  t_i  ] \prod_{j=1}^{i-1} \prob [ \tau_j > t_i ] \Big).
\end{align}
Note that $\expt[\tau_i | x(\tau_i)=\zeta \; \&\; \tau_i \le t_i ]
\er_i = \expt[\tau_i \bs 1 (x(\tau_i)=\zeta) | \tau_i \le t_i ]$ and
$\prob [ \tau_i \le t_i ] \prod_{j=1}^{i-1} \prob [ \tau_j > t_i ] =
\prob [ t_{i-1} < \tau \le t_i ]$.  Similar to
\eqref{eq:multistage-DT}, these equations show that the conditional
decision times are weighted sums of the expected conditional decision
times for each stage, with the weights being the conditional
probability of the decision in each stage.

\item The FPT cumulative distribution functions conditioned on hitting  upper/lower threshold are  
 \begin{align}
\prob [\tau \le t| x(\tau)=\zeta]  = \; & \frac{1}{\subscr{P}{ms}^+}\bigg(\prob[\tau_k \le t \; \& \; x(\tau_k) =\zeta]\prod_{j=1}^{k-1} \prob [ \tau_j > t_j ]  \nonumber\\ \label{eq:multistage-fpt-correct}
&\;\; \qquad \qquad + \sum_{i=1}^{k-1} \prob [\tau_i \le t_i \; \& \; x(\tau_i)=\zeta] \prod_{j=1}^{i-1} \prob [ \tau_j > t_j ]
\bigg)\\
\prob [\tau \le t| x(\tau)=-\zeta]  = \; & \frac{1}{\subscr{P}{ms}^-}\bigg(\prob(\tau_k \le t \; \& \; x(\tau_k) =-\zeta) \prod_{j=1}^{k-1} \prob [ \tau_j > t_j ]\nonumber \\ \label{eq:multistage-fpt-incorrect} 
& \qquad \qquad + \sum_{i=1}^{k-1} \prob [\tau_i \le t_i \; \& \; x(\tau_i)=-\zeta] \prod_{j=1}^{i-1} \prob [ \tau_j > t_j ]\bigg).
\end{align}
Note that $\prod_{j=1}^{i-1} \prob [ \tau_j > t_j ] = \prob [ \tau > t_{i-1} ]$ and $\prob [\tau_i \le t_i \; \& \; x(\tau_i)= \zeta]\prob [ \tau > t_{i-1} ] = \prob [t_{i-1} < \tau \le t_i \; \& \; x(\tau)=\zeta]$.
\end{enumerate}

\section{{\sfr Time-varying thresholds for the multistage process}}
\label{sec:thresholds}
\begin{figure}
\centering
\includegraphics[width=0.5\textwidth]{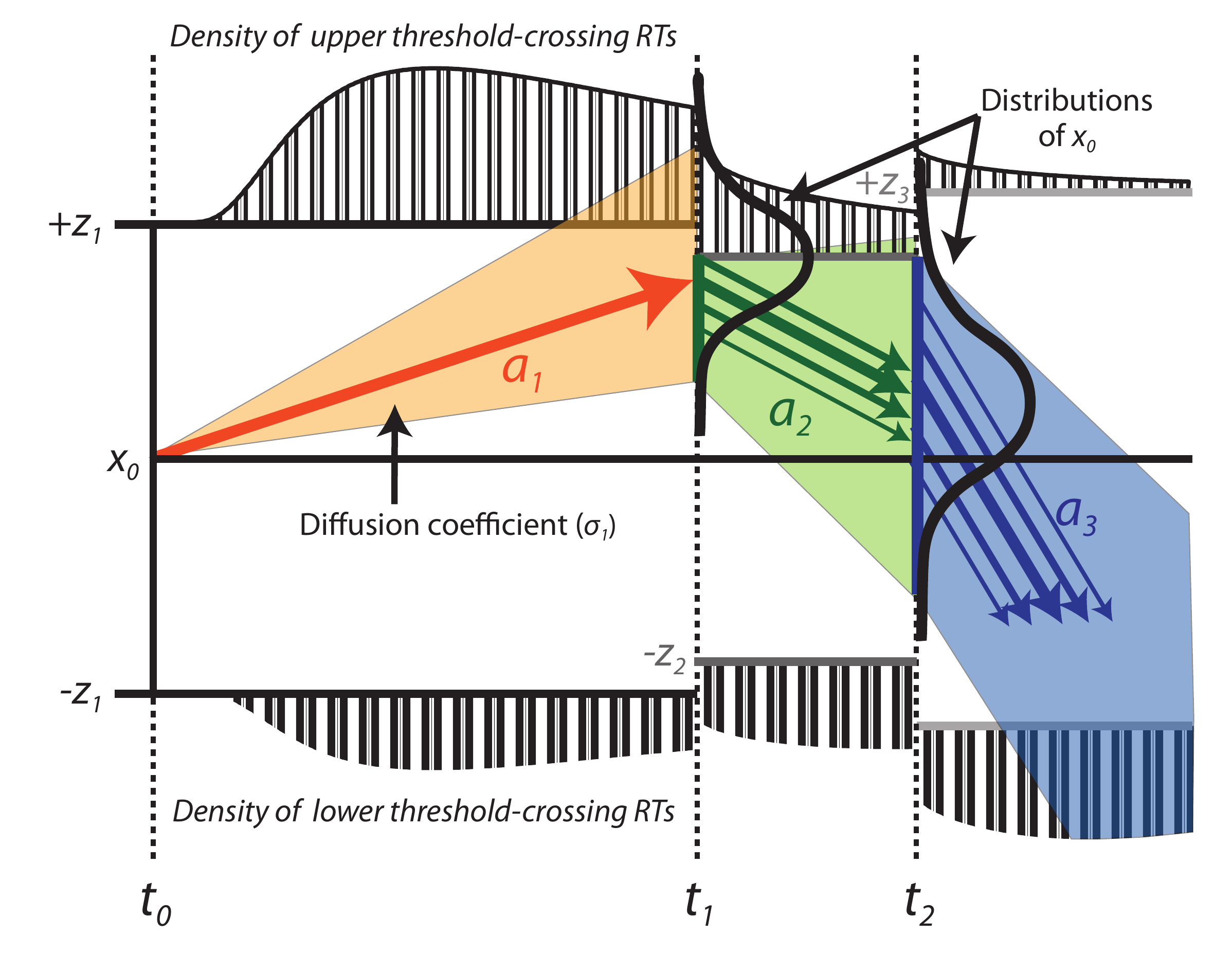}
\caption{Illustration of the key ideas for computation of FPT properties of {\sfr the multistage model with piecewise constant drift and thresholds}. The distribution of evidence $x(t_i)$ conditioned on no decision until {\sfr $t_i$ serves as} the distribution of initial condition for the $i$-th stage {\sfr process}. If threshold $\zeta_{i+1} <\zeta_i$, then probability {\sfr of} instantaneous decision at time $t_i$ is computed as the probability of $x(t_i)$ conditioned on no decision until $t_i$ not belonging to the set $(-\zeta_i, \zeta_i)$. If threshold $\zeta_{i+1} >\zeta_i$, then there is no instantaneous decision and only the support of  $x(t_i)$ conditioned on no decision until $t_i$ is increased to $(-\zeta_{i+1}, \zeta_{i+1})$. \label{fig:msddm-schematic}}
\end{figure}

The results in \S\ref{sec:performance-metrics} were obtained under the assumption that the thresholds are constant throughout each stage. 
Now suppose that the thresholds for the $i$-th {\sfr stage are} $\pm \zeta_i$, i.e., piecewise constant thresholds.
{\sfr If the upper thresholds decrease at time $t_i$ (i.e. $\zeta_{i+1} < \zeta_i$) and $x(t_i)$ is in the interval $(\zeta_{i+1},\zeta_i)$, then the path is absorbed by the upper boundary, and the probability of this instantaneous absorption is calculated by integrating \eqref{eq:dist-deadline} from $\zeta_{i+1}$ to $\zeta_i$.}
Likewise, the probability of instantaneous absorption into the lower threshold at $t_i$ is determined by integrating \eqref{eq:dist-deadline} from $-\zeta_i$ to $- \zeta_{i+1}$.
The density of $X_i$ is then found by truncating the support of the density in \eqref{eq:dist-deadline} to $(-\zeta_{i+1}, \zeta_{i+1})$ and normalizing the truncated density. 
In the cases where the upper threshold in the $(i+1)$-th stage is larger than the upper threshold in the $i$-th stage, i.e., $\zeta_{i+1} > \zeta_i$, there is no instantaneous absorption, and the density of $X_i$ is found by extending the density in \eqref{eq:dist-deadline}, assigning zero density to the previously undefined support (see Figure~\ref{fig:msddm-schematic}).
In all cases, the new, updated $X_i$ may be used for computations dealing with the $(i+1)$-th stage of the {\sfr multistage model}.
Codes implementing all of the formulas through \S\ref{sec:performance-metrics} with time-varying thresholds may be found at \texttt{https://github.com/sffeng/multistage}.  In \ref{app:ou}, we describe how these ideas extend to a time varying Ornstein-Uhlenbeck (O-U) model.

\section{Numerical examples} \label{sec:numerics}
In this section we apply our calculations from
\S\ref{sec:performance-metrics} and \S\ref{sec:thresholds} to a
variety of numerical experiments.  In doing so, we compare the
theoretical predictions obtained from the analysis in this paper with
the numerical values obtained through Monte-Carlo simulations, thereby
numerically verifying our derivations above.  We also provide examples
illustrating time pressure or changes in attention over the course of
a decision process, and demonstrate how our work can help to find the
optimal speed-accuracy trade-off by maximizing reward rate (or any
other function of mean first passage time, threshold-hitting
probability, and reward). Unless otherwise noted, Monte Carlo
simulations were obtained using $1000$ runs; relatively few runs are
used so that curves are visually distinguishable.  {\sfr Stochastic
  simulations use the} Euler-Maruyama method with time step size
$10^{-3}$.  All of the above calculations have been implemented in
MATLAB, and all codes used to produce the figures in this section may
be found at \texttt{https://github.com/sffeng/multistage}.

\subsection{A Four-stage {\sfr process}}\label{sec:num-exp}
\begin{figure}
\centering 
\subfigure[Unconditional and Conditional First Passage Time Distribution]{\label{fig:fpt-dist}
\includegraphics[width=0.7\linewidth]{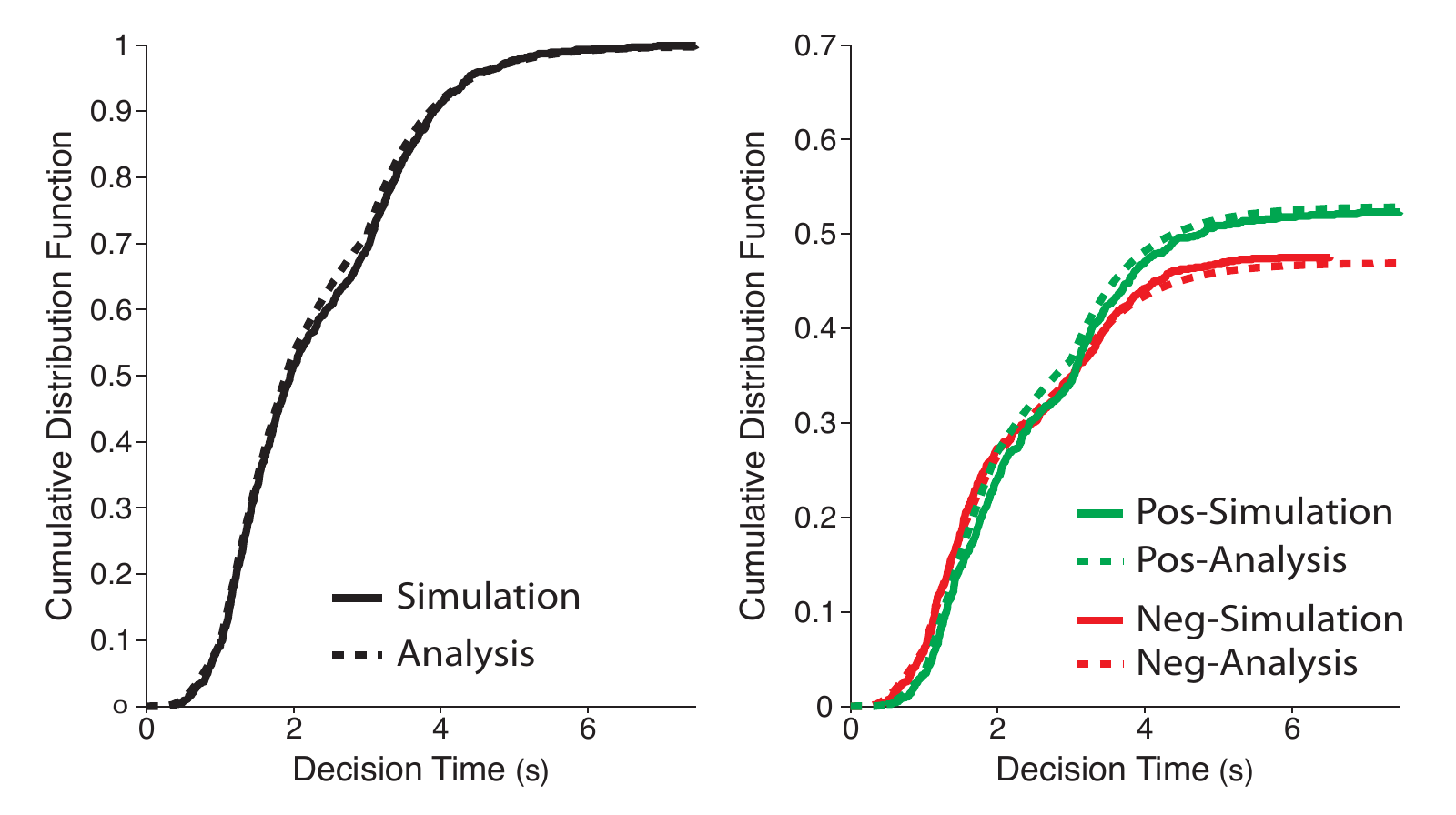}}
\subfigure[Unconditional Mean Decision Time, Lower Hitting Probability, and Conditional Mean Decision Time]{\label{fig:decision-time}
\includegraphics[width=0.7\linewidth]{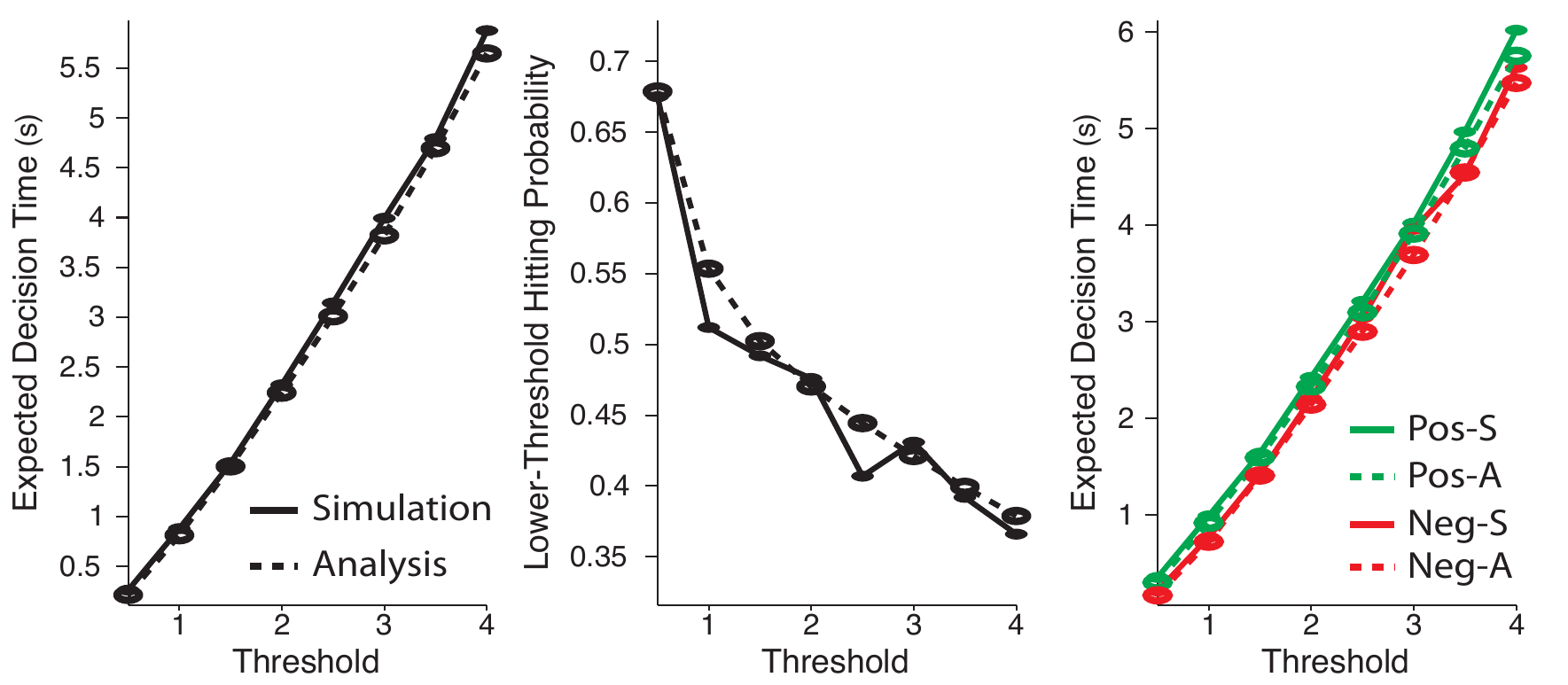}
}
\caption{FPT calculations for a four-stage {\sfr process} with drift rates $(\mu_1, \mu_2, \mu_3, \mu_4)= (0.1, 0.2, 0.05, 0.3)$, diffusion rates $(\sigma_1, \sigma_2, \sigma_3, \sigma_4) = (1, 1.5, 1.25, 2)$, stage initiation times $(t_0, t_1, t_2, t_3) =(0, 1 ,2, 3)$, and initial condition $x_0=-0.2$. The FPT distribution is computed for threshold $\zeta=2$. \label{fig:four-stage-ddm}}
\end{figure}

Consider a four stage {\sfr process} with drift and diffusion rates given by $(\mu_1, \mu_2, \mu_3, \mu_4)= (0.1, 0.2, 0.05, 0.3)$ and $(\sigma_1, \sigma_2, \sigma_3, \sigma_4) = (1, 1.5, 1.25, 2)$, respectively, with $(t_0, t_1, t_2, t_3) =(0, 1 ,2, 3)$ and initial condition $x_0=-0.2$. 
The cumulative distribution function (CDF) of the unconditional and conditional decision time for $\zeta=2$ obtained using the above analytic expressions (solid lines) and Monte-Carlo simulations (dotted lines) is shown in Figure~\ref{fig:fpt-dist}.
Similarly, the unconditional mean decision time, the lower hitting probability, and the mean conditional decision times (for upper/lower responses) are shown in Figure~\ref{fig:decision-time} as a function of threshold $\zeta$. {Note that the analytic expressions match closely with quantities computed using Monte-Carlo simulations. Also, notice that the CDF almost  looks like a double-sigmoidal function (it starts to saturate around $0.6$ before picking up and eventually saturating at $1$) due to the drop in drift rate from $0.2$ to $0.05$.}

\subsection{FPT Distribution for a {\sfr process} with Alternating Drift}\label{sec:num_kraj}
\begin{figure}
\centering
\includegraphics[width=0.7\textwidth]{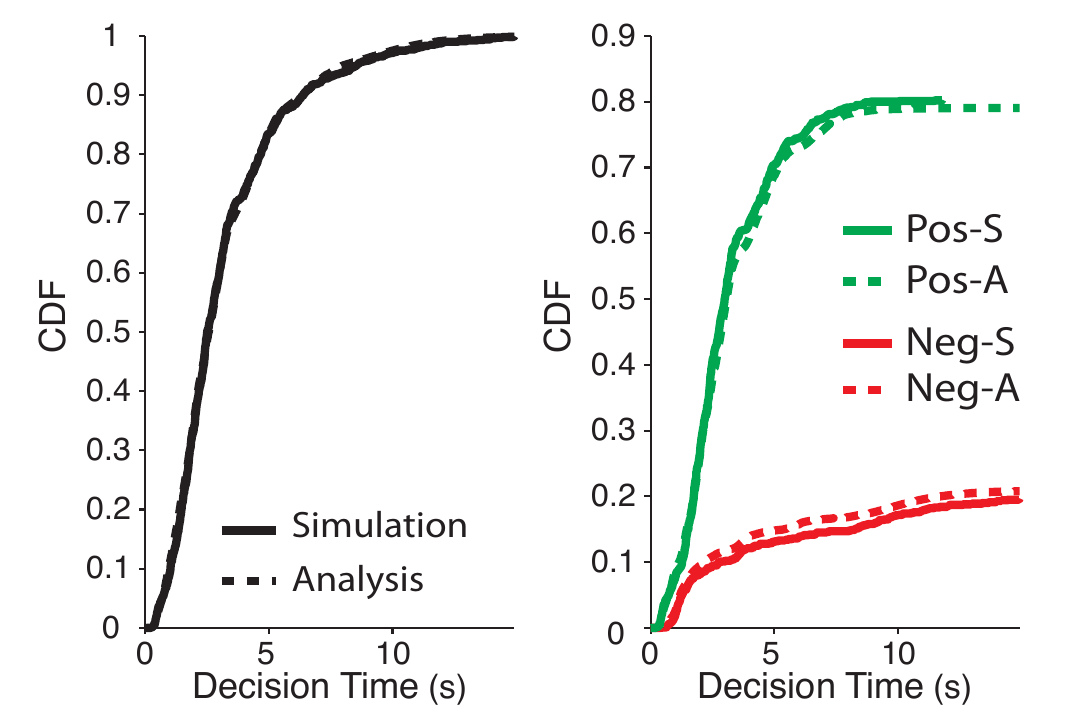}
\caption{Unconditional and conditional FPT distributions for a
  $30$-stage {\sfr model} with alternating drift rate. The drift rates are
  $(\mu_1, \mu_2,\mu_3, \ldots)=(1, -0.75, 1, \ldots)$, diffusion rate at
  each stage is unity, the threshold $\zeta=2$, and stage initiation times
  are {\sfr equally spaced throughout} the interval $(0, 10)$.
\label{fig:Krajbich}}
\end{figure}

In this example, the sign of the drift rate changes from stage to
stage.  This may be used to describe situations in which evidence
accumulation changes dynamically with the decision-maker's focus of
attention. For instance, \cite{IK-CA-AR:10} have shown that the
process of weighing two value-based options (e.g., foods) can be
modeled with a process in which drift rates vary based on the option
being attended at any given moment. We consider such a case using a
$30$-stage {\sfr model} in which the drift rates $1$ and $-0.75$
alternate (i.e., $\mu_1=1$, $\mu_2=-0.75$, $\mu_3=1$, $\ldots$) to
capture a situation in which the decision maker's attention alternates
between two options, one of which has greater perceived value (higher
drift rate) than the other. Let $t_0=0$ and the remaining $29$ stage
initiation times be a fixed realization of $29$ uniformly sampled
points between $0$ and $10$.  Assume $x_0 =0$, $\zeta=2$, and let the
diffusion rate be stationary and equal to unity ($\sigma_i = 1$). The
unconditional and conditional FPT distributions in this scenario
obtained using both the analytic expressions (solid lines) and
Monte-Carlo simulations (dotted lines) are shown
in~Figure~\ref{fig:Krajbich}. {Note that the analytic expressions
  match closely with quantities computed using Monte-Carlo
  simulations.}

\subsection{FPT Distribution for a {\sfr model with gradually time-varying} drift}\label{sec:num-white}

Changes in evidence accumulation may occur gradually over time. For
instance, \cite{CNW-RR-JFS:11} proposed a ``shrinking spotlight" model
of the Eriksen Flanker Task, a task in which participants responding
to the direction of a central arrow are influenced by the direction of
arrows in the periphery (see also \citealp{DSS-HP-JC:90, SL-AJY-PH:09,
  MS-CW-AM-BB:15}).  According to these models, evidence accumulation
is initially influenced by all of the arrows (central as well as
flankers, which may drive an incorrect response, modeled here as as a
lower threshold response) but as the attentional spotlight narrows the
drift rate is gradually more influenced by the central arrow alone.
{\sfr The multistage model, in spite of having discontinuous changes
  in parameters, can still be used to approximate a model with
  gradually time-varying parameters.}

\begin{figure}
\centering
\includegraphics[width=0.7\textwidth]{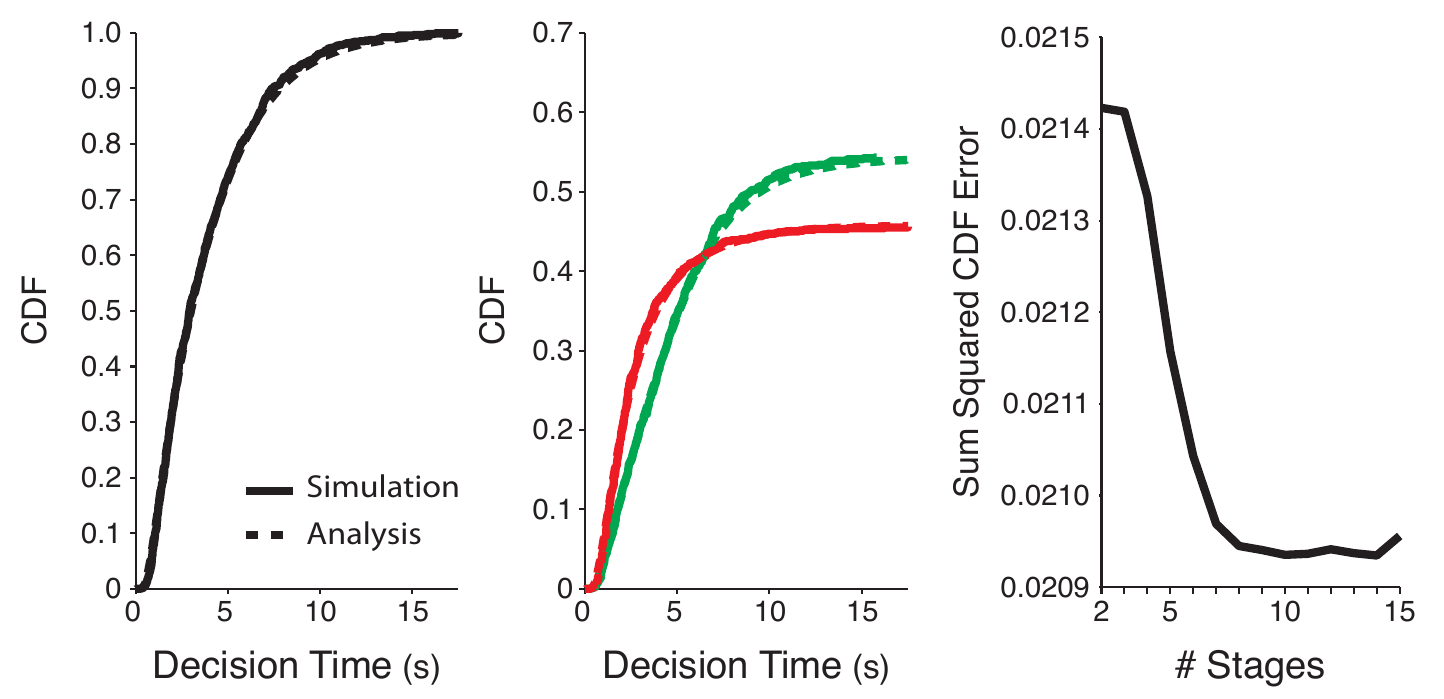}
\caption{Unconditional (left) and conditional (middle) FPT distributions for a $20$-stage {\sfr process} with gradually increasing drift rate. The drift rate for the $i$-th stage is  $\mu_i = -0.2 + 0.0263(i-1)$, diffusion rate at each stage is unity, the threshold $\zeta=2$, and stage initiation times are {\sfr equally spaced throughout} the interval $[0, 5]$. Right: The total squared error between the simulated and analytic CDFs (left) decreases with increasing discrete stages of the {\sfr model}. 10,000 simulations were used for each approximation.
\label{fig:gradual-time-varying}}
\end{figure}

As a demonstration, we use a $20$-stage {\sfr model} as an
approximation to {\sfr a model with continuously} time-varying drift
rate. Assume $\sigma_i = 1$, $x_0=0$, $\zeta=2$, and let the stage
initiation times $t_0, t_1, \ldots, t_{19}$ be {\sfr equally spaced
  throughout} the interval $[0,5]$. Furthermore, suppose the drift
rate during the $i$-th stage is $-0.2 + 0.0263(i-1)$. The
unconditional and conditional FPT distributions for such a $20$-stage
{\sfr process} obtained using the analytic expressions (solid lines)
and using Monte-Carlo simulations (dotted lines) are shown in
Figure~\ref{fig:gradual-time-varying}. {Note that the analytic
  expressions match closely with quantities computed using Monte-Carlo
  simulations.} We also show the error due to piecewise constant
approximation of the drift rate as a function of number of stages (Figure~\ref{fig:gradual-time-varying}, right). It
should be noted that even for $5$ stages the approximation error is very
small.

\subsection{Collapsing Thresholds}\label{sec:num-collapse}
\begin{figure}
\centering
\includegraphics[width=0.7\textwidth]{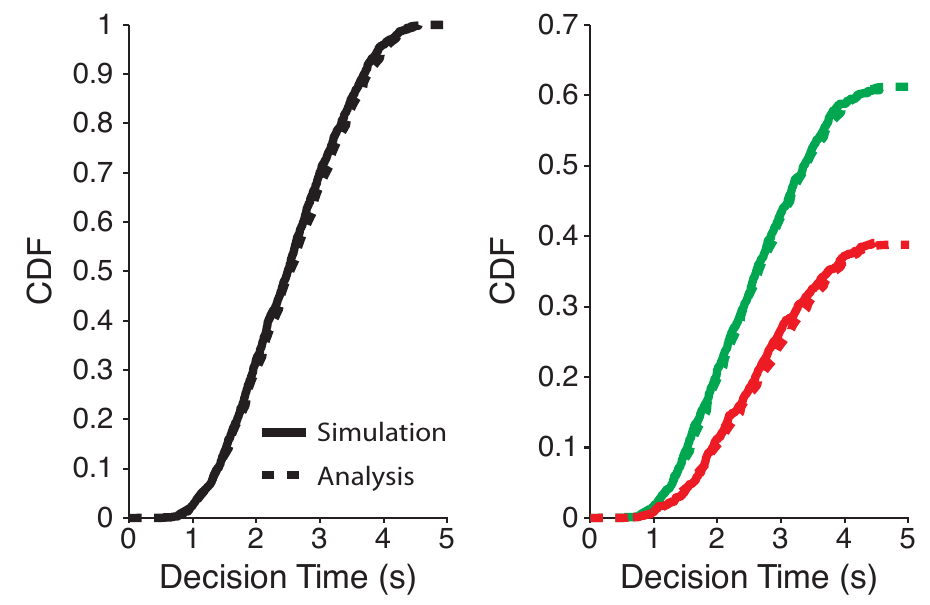}
\caption{Unconditional and conditional FPT distributions for a $20$-stage {\sfr process} with collapsing thresholds. The drift rate and diffusion rate at each stage are $0.15$ and $1$, respectively.  The stage initiation times are {\sfr equally spaced throughout} the interval $[0, 5]$, and the threshold at $i$-th stage is $\zeta_i= 3-\frac{i-1}{19}$, for $i = \until{19}$.
\label{fig:collapsing-threshold}}
\end{figure}

One may be interested in modeling a decision process in which
thresholds are dynamic rather than {\sfr constant} across stages. This
can be used to describe discrete changes in choice strategy, or a
continuous change in thresholds over time; the latter approach has
been successful at describing behavior under conditions that either
involve an explicit response deadline (e.g.,~\citealp{MM-JM-etal:10,
  PF-AJY:08}) or where there is an implicit opportunity cost for
longer time spent accumulating evidence~\citep{JD-RMB-etal:12}. {\sfr
  Recently, \citep{Voskuilen2016} used analytic methods to model
  collapsing boundaries in order to compare fixed boundaries against
  collapsing boundaries in diffusion models.

Here, we } model such a situation using a $20$-stage {\sfr process}, as an
approximation to a diffusion model with continuously collapsing
thresholds ,i.e., $\zeta \downarrow 0$ with time, the drift rate and the
diffusion rate are constant and equal to $0.15$ and $1$, respectively,
$x_0=0$, and stage initiation times $t_0,t_1,\ldots,t_{19}$ are
{\sfr equally spaced throughout} the interval $[0,5]$. The threshold in the
$i$-th stage is $\zeta_i= 3- \frac{1}{19}(i-1)$. The unconditional and
conditional FPT distributions for such a $20$-stage {\sfr process} obtained
using the analytic expressions (solid lines) and using Monte-Carlo
simulations (dotted lines) are shown in
Figure~\ref{fig:collapsing-threshold}. {Note that the analytic
  expressions match closely with quantities computed using Monte-Carlo
  simulations.}

\subsection{Optimizing the Speed-Accuracy Trade-off in a Two-stage {\sfr Model}}\label{sec:speed-accuracy}
Human decision-making in many two alternative forced choice signal detection tasks has been successfully captured by the {\sfr single stage model}. In such tasks, hitting the upper/lower boundary is interpreted as a correct/incorrect response.  The accuracy of a decision can then be determined by the sign of the drift rate -- if $\mu$ is positive, participants are said to be more accurate the more likely they are to hit the upper threshold and more error-prone the more likely they are to hit the lower threshold.
One then assumes without loss of generality that the drift rate is positive, in which case the lower hitting probability is called the \emph{error rate} and the upper hitting probability is called the \emph{accuracy}.
Here, the choice of threshold dictates the speed-accuracy trade-off, i.e., the trade-off between a fast decision and an accurate decision. 
In the previous examples, 
the thresholds have been known and we have characterized the associated error rate and first passage time properties. 
These can be used to define a joint function of speed and  accuracy that may dictate how humans/animals choose to set and adjust their threshold. In particular, it has been proposed~\citep{Bogacz2006} that human subjects choose a threshold that maximizes reward rate (RR), defined as
\begin{equation}\label{eq:rr}
{\rm RR} = \frac{1 - \er}{ \expt[\tau] + \subscr{T}{nd}},
\end{equation}
where $ \subscr{T}{nd}$ is the sensory and motor processing time (non-decision time) and $\er$ and $\expt[\tau]$ are computed using the expressions derived in \S\ref{sec:performance-metrics}.

\begin{figure}[ht!]
\centering
\includegraphics[width=0.495\textwidth]{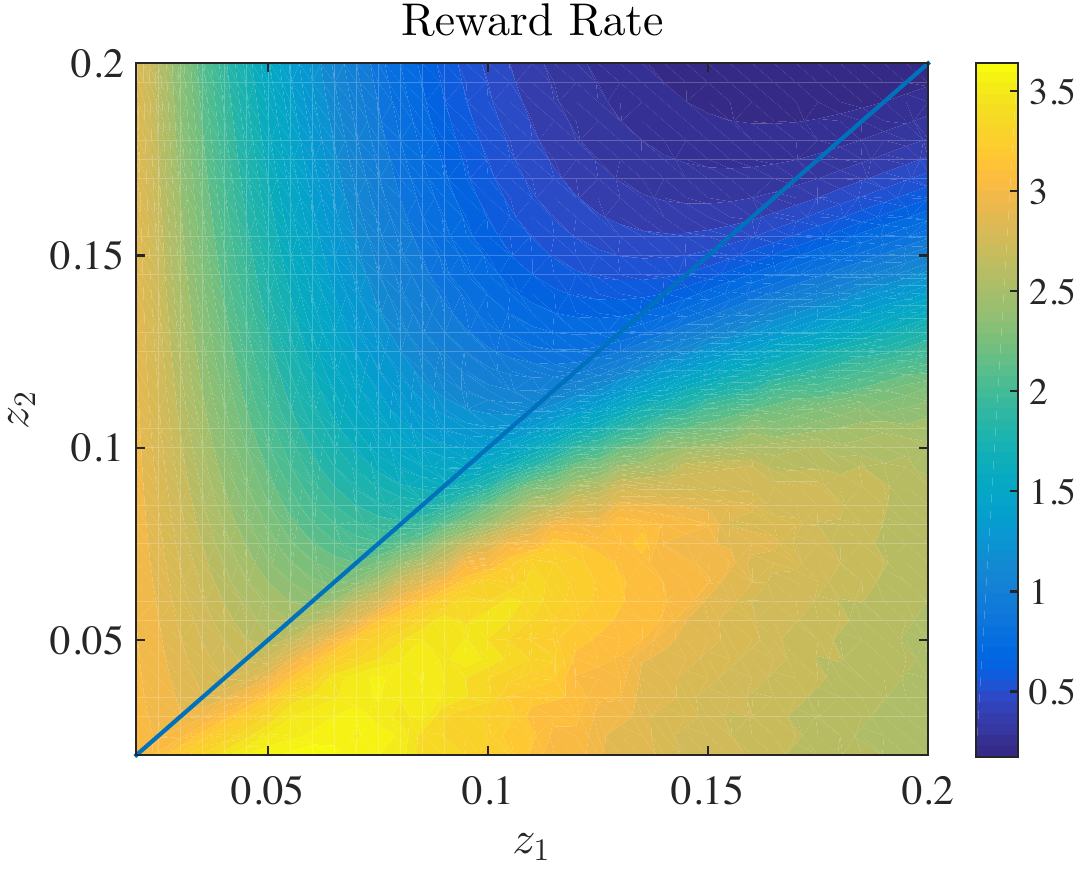}
\caption{Reward rate as a function of the thresholds
 for the two-stage {\sfr model}. The {\sfr parameters are} $\mu_1=0.5$, $\mu_2=-0.1$, $\sigma_1=\sigma_2=0.1$, $x_0=\zeta_1/2$, and $t_1=0.1$. The non-decision time $\subscr{T}{nd}=0.3$.
Note that the maximum reward rate is achieved for $\zeta_1 \ne \zeta_2$. \label{fig:rr-surf}}
\end{figure}

The reward rate for a two-stage {\sfr process} is shown in Figure~\ref{fig:rr-surf}. For the set of parameters in Figure~\ref{fig:rr-surf}, reward rate is maximal at approximately $(\zeta_1, \zeta_2)=(0.06, 0.01)$. Thus, the maximizing reward rate in this setting interestingly requires that the threshold across stages be different. 

\begin{figure}[ht!]
\centering
\subfigure[Reward rate for single- and two-stage {\sfr models}]{\includegraphics[width=0.475\textwidth]{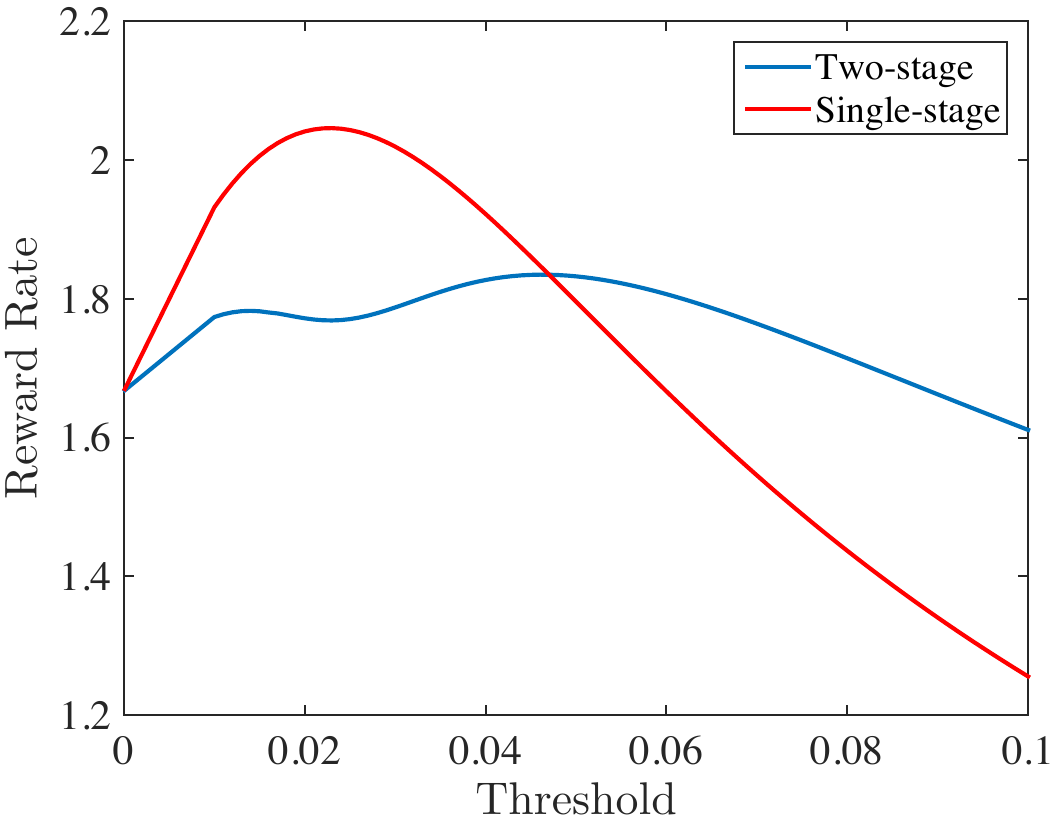}\label{fig:rr-2-ddm}}
\subfigure[Reward rate for two-stage {\sfr model}]{\includegraphics[width=0.475\textwidth, height=0.4\textwidth]{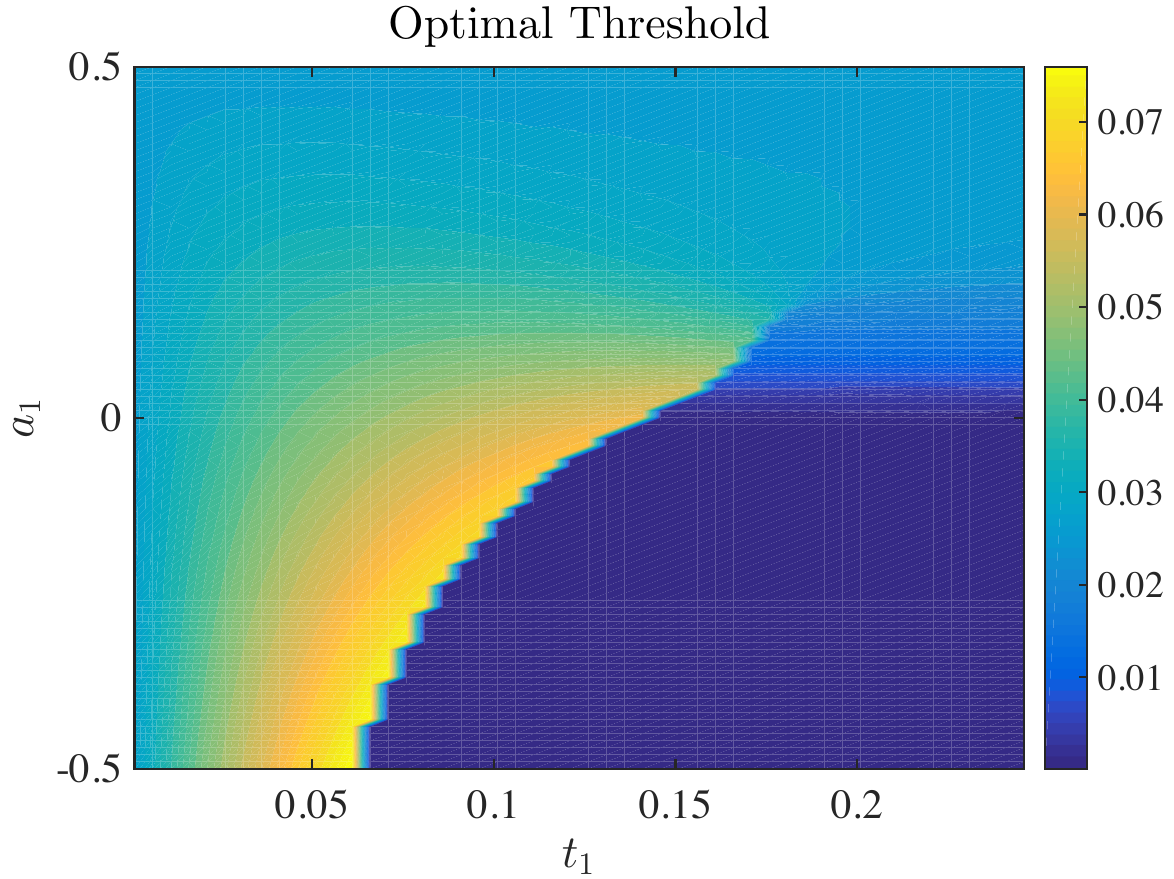} \label{fig:opt-thresh-rr}}
\caption{
(a) For the single-stage {\sfr model} $\mu_1=0.2, \sigma_1=0.1$ and $x_0=0$, while for the two-stage {\sfr model} $\mu_1 =0.1$, $\mu_2=0.5$, $\sigma_1=\sigma_2=0.1$, $x_0=0$ and $t_1= 0.15$. The non-decision time $\subscr{T}{nd}= 0.3$. The reward rate for the single-stage {\sfr model} is a unimodal function and achieves a unique local maximum, while the reward rate for the two-stage {\sfr model} has two local maxima. 
(b) Optimal threshold for two-stage {\sfr model} obtained by maximizing reward rate. The left panel shows the variation of the optimal threshold as a function of $t_1$ and $\mu_1$. The other parameters are  $\mu_2=0.5$, $\sigma_1=\sigma_2=0.1$, and $x_0=0$. The regions of the contour plot associated with $t_1=0$ and $\mu_1=0.5$ correspond to the single-stage {\sfr model}.}
\end{figure}

Setting the threshold to be constant across stages ($\zeta_1=\zeta_2$), we can compare how reward rate changes with this constant threshold in a single-stage (traditional) {\sfr versus a two-stage mode}. As shown in Figure~\ref{fig:rr-2-ddm}, we find that this reward rate function is unimodal {\sfr in a single-stage model (as previously observed) whereas it is bimodal in a two-stage model}. Figure~\ref{fig:opt-thresh-rr} explores this parameter space in greater depth and shows that the curvature of reward rate (and in particular the relative height of its first and second modes) vary as a function of the length and drift rate of the first stage (for example). As a result, this analysis reveals a discontinuous jump in optimal threshold as these parameters vary. Whether individuals are sensitive to these discontinuities when setting thresholds for a multistage decision-making task deserves further exploration.

\section{Discussion}
\label{sec:discussion}
\label{sec:conclusion}

In this paper we analyze the first passage time properties of a Wiener
process between two absorbing boundaries with piecewise constant
(time-dependent) parameters, which we call a multistage model or
multistage process. Our main theoretical results, collected in \S
\ref{sec:performance-metrics}, add to previous work on analyzing
time-dependent random walk models in psychology and
neuroscience. Broadly speaking, these can be split into three
approaches. One approach is the integral equation approach introduced
and developed in \cite{Smith1995, Smith2000, Smith2009}. Another
approach is the matrix based Markov Chain approximation which has been
applied to a wide variety of multi-attribute choice settings
\citep{Diederich2003, Diederich2016, Diederich2014}. A third approach
analyzes the backward partial differential equation associated with
the multistage process \citep{Ratcliff1980,Heath1992}.  The results of
\S \ref{sec:performance-metrics} align most closely with the work of
the third approach, as we also directly analyze the multistage
stochastic process, albeit with different techniques.  Whereas
\cite{Ratcliff1980} and \cite{Heath1992} analyzed a multistage process
by solving the Kolmogorov partial differential equation, we employ
martingale theory (e.g. the optional sampling theorem) in order to
obtain analytic results which extend those of previous studies.  In
doing so, our work also builds on martingale-based analyses described
by \cite{smith1990note} for a single stage model.

The modeler utilizing time-dependent random walk models in decision
making should be aware of all of the above approaches, as one approach
may demand additional approximations compared to another, due to
differences in how they discretize temporal dynamics. For example,
when modeling experiments with continuously (gradually) changing
stimuli (e.g., \cite{CNW-RR-JFS:11, zhang2014time}) one may find the
integral equation approach more natural, given the continuity and
smoothness assumptions built in to the techniques. Diederich and
Oswald's model also allows for continuous changes in boundaries. In
contrast, our methods must approximate the underlying gradually
changing drift with a piecewise constant function, thereby inducing
some additional error in the calculations of first passage times.  If,
however, the application lends itself well to discrete changes in
drift rate or threshold (e.g., \cite{IK-CA-AR:10}), then both the
matrix based methods and methods discussed in this paper may be more
natural, as they explicitly consider such discrete changes in the
underlying calculations.  Ornstein-Uhlenbeck processes lend themselves
more naturally to the matrix based approach of \cite{Diederich2003},
compared to our analysis in \ref{app:ou}, since our analysis requires
a change in time scale that the matrix method does not.  In the
specific case of a multistage model \eqref{eq:multistage-ddm}, our
work provides semi-analytic formulas that can be easily computed or
studied further (see below). In general, the modeler should be ready
to employ the most suitable approach given their situation.

The reader should also be aware of the software tools available for
each approach, as there already exist several good non-martingale
based software packages for computing FPT statistics. One package
implementing the integral equation approach is that of~\cite{JD:14},
which computes first passage time densities using the stable numerical
approximations developed in~\cite{Smith2000}. More recently, highly
optimized codes for a broad class of diffusion models have been
developed by \cite{Verdonck2015}, with implementations on both CPUs
and GPUs.  Compared to other available codes, \cite{Diederich2003}'s
matrix approach is substantially simpler and more elegant to implement
-- one can compute desired choice probabilities and mean decision
times in less than a few dozen lines of MATLAB code.  For
practitioners wishing to write all of their own model code from
scratch, this may be a considerable advantage. The MATLAB code
released with this
report\footnote{\texttt{https://github.com/sffeng/multistage}} provides
implementations of the results from \S \ref{sec:performance-metrics}
and allows one to reproduce all of the figures from \S
\ref{sec:numerics}.  Unlike the work of \cite{Verdonck2015}, our work
is not immediately focused on developing a rapid numerical tool for
simulation, but rather on introducing martingale theory as a useful
mathematical tool for analyzing multistage decision models.  Thus, the
codes released with this report are not intended to compete with the
efficiency of the aforementioned codes, which have been highly
optimized and tuned for throughput, but instead to demonstrate the
simplicity and effectiveness of our analysis.  These considerations
notwithstanding, our results do suggest promising avenues for future
numerical work.  Particularly relevant is work by \cite{Navarro2009,
  Blurton2012, Gondan2014} who developed efficient numerical schemes
for evaluating the relevant infinite sums involved in FPT
calculations.  Similar methods could be applied to results in
\S\ref{subsec:i-ddm}~and~\S\ref{subsec:multistage-ddm} to develop
efficient multistage codes, which could in turn contribute to the
growing collection of numerical tools available for practitioners
using diffusion models to study decision making.

In \S\ref{sec:review} we re-derived classical mean decision times and
choice probabilities equivalent to results first derived for a
discrete time random walk model using the Wald identity
\citep{Laming1968, Link1975, Link1975a, smith1990note}, which itself
is a corollary of the optional sampling theorem.  \cite{smith1990note}
notes that results concerning the discrete time random walk model with
Gaussian increments and the continuous time single stage Wiener
diffusion model of \S\ref{sec:review} should be equivalent, which is
indeed the case. Furthermore, the moment generating function derived
in \S\ref{subsec:app-ddm} is identical to that obtained by
\cite[p. 9]{smith1990note} via the Wald identity.  Our aim in
presenting \S\ref{sec:review} was to introduce the martingale analysis
used in \S\ref{sec:performance-metrics} by first demonstrating its
utility in re-deriving these classical results. We hope these
calculations provide an intuition for the computations in
\S\ref{sec:performance-metrics} in a less technical setting.  For
reviews of the classical results on random walks in psychological
decision making, see \cite[pp. 300-301]{Townsend1983} and
\cite[pp. 328-334]{Luce1986}.

One particularly useful aspect of equations
\eqref{eq:error-rate-deadline} through
\eqref{eq:multistage-fpt-incorrect} is that they demonstrate how
various other measures of performance, such as the lower-threshold
hitting probability during each stage, evolve as the underlying
dynamics change.  Using these, one may efficiently compute a variety
of behavioral measures of performance without resorting to first
computing the FPT densities. Our results may also serve as a starting
point for further analysis of more complicated stochastic decision
models.  In appendix D, we show how our results apply to
Ornstein-Uhlenbeck processes, which approximate leaky integration over
the course of evidence accumulation, e.g., the Leaky Competing
Accumulator model [LCA]~\citep{Usher2001}. Given that the LCA itself
can in certain cases approximate a reduced form of more complex and
biologically plausible models of interactions across neuronal
populations (e.g., \citealp{XJW:02,KFW-XJW:06,RB:07}), {\sfr our work
  may help analysts better understand time-varying dynamics within and
  across neural networks, and how such dynamics relate to complex
  cognitive phenomena}.

Beyond its analytic and numerical utility, our uses of martingale
theory, and the optional sampling theorem in particular, provide some
theoretical insights into random walk decision models. For example, a
single stage model with zero drift rate (i.e. Wiener process) and
lower/upper thresholds at 0 and 1 necessarily has the qualitative
property that the probability of hitting the upper threshold equals
the initial (nonrandom) position $x_0$. Although these results are
known in the probability literature~\citep{Doob1953}, they provide
practitioners and experimentalists a unique way of planning and
analyzing experiments based on decisions thought to evolve according
to diffusion processes. In summary, we hope that the tools described
in this article and the results of \S \ref{sec:performance-metrics}
will encourage computational and mathematical analyses of decision
models involving time-dependent parameters.

\section{Acknowledgements}
We also thank the editor, Philip Smith, and the referees for helpful
comments, especially regarding this article's exposition and
discussion.  We thank Ryan Webb for the reference to~\cite{Smith2000}
and the associated code.  We thank Phil Holmes and Patrick Simen for
helpful comments and discussions.  This work was jointly supported by
the C.V. Starr Foundation (AS), the Insley-Blair Pyne Fund, ONR grant
N00014-14-1-0635 and ARO grant W911NF-14-1-0431(VS and NEL).

\section{Bibliography}

{\footnotesize
\bibliographystyle{elsarticle-harv.bst}
\bibliography{msddm}

\begin{thebibliography}{64}
\expandafter\ifx\csname natexlab\endcsname\relax\def\natexlab#1{#1}\fi
\expandafter\ifx\csname url\endcsname\relax
  \def\url#1{\texttt{#1}}\fi
\expandafter\ifx\csname urlprefix\endcsname\relax\def\urlprefix{URL }\fi

\bibitem[{Blurton et~al.(2012)Blurton, Kesselmeier, and Gondan}]{Blurton2012}
Blurton, S.~P., Kesselmeier, M., Gondan, M., 2012. Fast and accurate
  calculations for cumulative first-passage time distributions in {W}iener
  diffusion models. Journal of Mathematical Psychology 56~(6), 470--475.

\bibitem[{Bogacz(2007)}]{RB:07}
Bogacz, R., 2007. Optimal decision-making theories: linking neurobiology with
  behaviour. Trends in Cognitive Sciences 11~(3), 118--125.

\bibitem[{Bogacz et~al.(2006)Bogacz, Brown, Moehlis, Holmes, and
  Cohen}]{Bogacz2006}
Bogacz, R., Brown, E., Moehlis, J., Holmes, P.~J., Cohen, J.~D., 2006. The
  physics of optimal decision making: {A} formal analysis of models of
  performance in two-alternative forced-choice tasks. Psychological Review
  113~(4), 700--765.

\bibitem[{Borodin and Salminen(2002)}]{Borodin2002}
Borodin, A.~N., Salminen, P., 2002. Handbook of {B}rownian Motion: Facts and
  Formulae. Springer.

\bibitem[{Brunton et~al.(2013)Brunton, Botvinick, and Brody}]{BWB-MMB-CDB:13}
Brunton, B.~W., Botvinick, M.~M., Brody, C.~D., 2013. Rats and humans can
  optimally accumulate evidence for decision-making. Science 340~(6128),
  95--98.

\bibitem[{Busemeyer and Diederich(2010)}]{Busemeyer2010}
Busemeyer, J.~R., Diederich, A., 2010. Cognitive modeling. Sage.

\bibitem[{Cisek et~al.(2009)Cisek, Puskas, and El-Murr}]{PC-GAP-SE:09}
Cisek, P., Puskas, G.~A., El-Murr, S., 2009. Decisions in changing conditions:
  The urgency-gating model. The Journal of Neuroscience 29~(37), 11560--11571.

\bibitem[{Cox and Miller(1965)}]{DRC-HDM:65}
Cox, D.~R., Miller, H.~D., 1965. The Theory of Stochastic Processes. Methuen \&
  Co. Ltd.

\bibitem[{Diederich and Busemeyer(2003)}]{Diederich2003}
Diederich, A., Busemeyer, J.~R., 2003. Simple matrix methods for analyzing
  diffusion models of choice probability, choice response time, and simple
  response time. Journal of Mathematical Psychology 47~(3), 304 -- 322.

\bibitem[{Diederich and Busemeyer(2006)}]{Diederich2006}
Diederich, A., Busemeyer, J.~R., 2006. Modeling the effects of payoff on
  response bias in a perceptual discrimination task: Bound-change,
  drift-rate-change, or two-stage-processing hypothesis. Perception {\&}
  Psychophysics 68~(2), 194--207.
\newline\urlprefix\url{http://dx.doi.org/10.3758/BF03193669}

\bibitem[{Diederich and Oswald(2014)}]{Diederich2014}
Diederich, A., Oswald, P., 2014. Sequential sampling model for multiattribute
  choice alternatives with random attention time and processing order.
  Frontiers in Human Neuroscience 8~(697), 1--13.

\bibitem[{Diederich and Oswald(2016)}]{Diederich2016}
Diederich, A., Oswald, P., 2016. Multi-stage sequential sampling models with
  finite or infinite time horizon and variable boundaries. Journal of
  Mathematical Psychology.\;In press.

\bibitem[{Doob(1953)}]{Doob1953}
Doob, J.-L., 1953. Stochastic Processes. John Wiley \& Sons, Inc., Chapman \&
  Hall.

\bibitem[{Douady(1999)}]{Douady1999}
Douady, R., 1999. Closed form formulas for exotic options and their lifetime
  distribution. International Journal of Theoretical and Applied Finance 2~(1),
  17--42.

\bibitem[{Drugowitsch(2014)}]{JD:14}
Drugowitsch, J., 2014. {C++} diffusion model toolset with {P}ython and {M}atlab
  interfaces. GitHub repository: \url{https://github.com/jdrugo/dm}, commit:
  {5729cd891b6ab37981ffacc02d04016870f0a998}.

\bibitem[{Drugowitsch et~al.(2012)Drugowitsch, Moreno-Bote, Churchland,
  Shadlen, and Pouget}]{JD-RMB-etal:12}
Drugowitsch, J., Moreno-Bote, R., Churchland, A.~K., Shadlen, M.~N., Pouget,
  A., 2012. The cost of accumulating evidence in perceptual decision making.
  The Journal of Neuroscience 32~(11), 3612--3628.

\bibitem[{Durrett(2010)}]{Durrett2010}
Durrett, R., 2010. Probability: Theory and Examples. Cambridge University
  Press.

\bibitem[{Farkas and Fulop(2001)}]{Farkas2001}
Farkas, Z., Fulop, T., 2001. One-dimensional drift-diffusion between two
  absorbing boundaries: Application to granular segregation. Journal of Physics
  A: Mathematical and General 34~(15), 3191--3198.

\bibitem[{Feller(1968)}]{Feller1968}
Feller, W., 1968. An Introduction to Probability Theory and its Applications.
  Vol.~1. John Wiley \& Sons.

\bibitem[{Feng et~al.(2009)Feng, Holmes, Rorie, and Newsome}]{SF-PH-AR-WTN:09}
Feng, S., Holmes, P., Rorie, A., Newsome, W.~T., 2009. Can monkeys choose
  optimally when faced with noisy stimuli and unequal rewards. PLoS
  Computational Biology 5~(2), e1000284.

\bibitem[{Frazier and Yu(2008)}]{PF-AJY:08}
Frazier, P., Yu, A.~J., 2008. Sequential hypothesis testing under stochastic
  deadlines. In: Platt, J., Koller, D., Singer, Y., Roweis, S. (Eds.), Advances
  in Neural Information Processing Systems 20. Curran Associates, Inc., pp.
  465--472.

\bibitem[{Gold and Shadlen(2001)}]{JIG-MNS:01}
Gold, J.~I., Shadlen, M.~N., 2001. Neural computations that underlie decisions
  about sensory stimuli. Trends in Cognitive Sciences 5~(1), 10--16.

\bibitem[{Gold and Shadlen(2007)}]{JIG-MNS:07}
Gold, J.~I., Shadlen, M.~N., 2007. The neural basis of decision making. Annual
  Review of Neuroscience 30~(1), 535--574.

\bibitem[{Gondan et~al.(2014)Gondan, Blurton, and Kesselmeier}]{Gondan2014}
Gondan, M., Blurton, S.~P., Kesselmeier, M., Jun. 2014. Even faster and even
  more accurate first-passage time densities and distributions for the {Wiener}
  diffusion model. Journal of Mathematical Psychology 60, 20--22.

\bibitem[{Heath(1992)}]{Heath1992}
Heath, R.~A., 1992. A general nonstationary diffusion model for two-choice
  decision-making. Mathematical Social Sciences 23~(3), 283--309.

\bibitem[{Horrocks and Thompson(2004)}]{Horrocks2004}
Horrocks, J., Thompson, M.~E., 2004. Modeling event times with multiple
  outcomes using the {Wiener} process with drift. Lifetime Data Analysis
  10~(1), 29--49.

\bibitem[{Hubner et~al.(2010)Hubner, Steinhauser, and Lehle}]{Hubner2010}
Hubner, R., Steinhauser, M., Lehle, C., 2010. A dual-stage two-phase model of
  selective attention. Psychological Review 117~(3), 759--784.

\bibitem[{Karatzas and Shreve(1998)}]{Karatzas1998}
Karatzas, I., Shreve, S., 1998. Brownian Motion and Stochastic Calculus, 2nd
  Edition. Vol. 113. Springer-Verlag New York.

\bibitem[{Krajbich et~al.(2010)Krajbich, Armel, and Rangel}]{IK-CA-AR:10}
Krajbich, I., Armel, C., Rangel, A., 2010. Visual fixations and the computation
  and comparison of value in simple choice. Nature Neuroscience 13~(10),
  1292--1298.

\bibitem[{Laming(1968)}]{Laming1968}
Laming, D. R.~J., 1968. Information theory of choice-reaction times. Academic
  Press.

\bibitem[{Lin(1998)}]{Lin1998}
Lin, X.~S., 1998. Double barrier hitting time distributions with applications
  to exotic options. Insurance: Mathematics and Economics 23~(1), 45--58.

\bibitem[{Link and Heath(1975)}]{Link1975}
Link, S., Heath, R., 1975. A sequential theory of psychological discrimination.
  Psychometrika 40~(1), 77--105.

\bibitem[{Link(1975)}]{Link1975a}
Link, S.~W., 1975. The relative judgment theory of two choice response time.
  Journal of Mathematical Psychology 12~(1), 114--135.

\bibitem[{Liu et~al.(2009)Liu, Yu, and Holmes}]{SL-AJY-PH:09}
Liu, S., Yu, A.~J., Holmes, P., 2009. Dynamical analysis of {B}ayesian
  inference models for the {E}riksen task. Neural Computation 21~(6),
  1520--1553.

\bibitem[{Luce(1986)}]{Luce1986}
Luce, R.~D., 1986. Response times: Their role in inferring elementary mental
  organization. No.~8. Oxford University Press.

\bibitem[{Milosavljevic et~al.(2010)Milosavljevic, Malmaud, Huth, Koch, and
  Rangel}]{MM-JM-etal:10}
Milosavljevic, M., Malmaud, J., Huth, A., Koch, C., Rangel, A., 2010. The drift
  diffusion model can account for the accuracy and reaction time of value-based
  choices under high and low time pressure. Judgment and Decision Making 5~(6),
  437--449.

\bibitem[{Mormann et~al.(2010)Mormann, Malmaud, Huth, Koch, and
  Rangel}]{mormann2010drift}
Mormann, M.~M., Malmaud, J., Huth, A., Koch, C., Rangel, A., 2010. The drift
  diffusion model can account for the accuracy and reaction time of value-based
  choices under high and low time pressure. Available at SSRN 1901533.

\bibitem[{Navarro and Fuss(2009)}]{Navarro2009}
Navarro, D.~J., Fuss, I.~G., 2009. Fast and accurate calculations for
  first-passage times in {W}iener diffusion models. Journal of Mathematical
  Psychology 53~(4), 222--230.

\bibitem[{Ratcliff(1978)}]{Rat78}
Ratcliff, R., 1978. A theory of memory retrieval. Psychological Review 85~(2),
  59--108.

\bibitem[{Ratcliff(1980)}]{Ratcliff1980}
Ratcliff, R., 1980. A note on modeling accumulation of information when the
  rate of accumulation changes over time. Journal of Mathematical Psychology
  21~(2), 178--184.

\bibitem[{Ratcliff and McKoon(2008)}]{Ratcliff2008}
Ratcliff, R., McKoon, G., 2008. The diffusion decision model: Theory and data
  for two-choice decision tasks. Neural Computation 20~(4), 873 -- 922.

\bibitem[{Ratcliff and Rouder(1998)}]{Ratcliff1998}
Ratcliff, R., Rouder, J.~N., 1998. Modeling response times for two-choice
  decisions. Psychological Science 9~(5), 347--356.

\bibitem[{Ratcliff and Smith(2004)}]{Ratcliff2004}
Ratcliff, R., Smith, P.~L., 2004. A comparison of sequential sampling models
  for two-choice reaction time. Psychological Review 111~(2), 333--367.

\bibitem[{Ratcliff et~al.(2016)Ratcliff, Smith, Brown, and
  McKoon}]{ratcliff2016diffusion}
Ratcliff, R., Smith, P.~L., Brown, S.~D., McKoon, G., 2016. Diffusion decision
  model: Current issues and history. Trends in Cognitive Sciences 20~(4),
  260--281.

\bibitem[{Revuz and Yor(1999)}]{Revuz1999}
Revuz, D., Yor, M., 1999. Continuous Martingales and Brownian Motion. Vol. 293.
  Springer Berlin Heidelberg.

\bibitem[{Servan-Schreiber et~al.(1990)Servan-Schreiber, Printz, and
  Cohen}]{DSS-HP-JC:90}
Servan-Schreiber, D., Printz, H., Cohen, J., 1990. A network model of
  catecholamine effects- gain, signal-to-noise ratio, and behavior. Science
  249~(4971), 892--895.

\bibitem[{Servant et~al.(2015)Servant, White, Montagnini, and
  Burle}]{MS-CW-AM-BB:15}
Servant, M., White, C., Montagnini, A., Burle, B., 2015. Using covert response
  activation to test latent assumptions of formal decision-making models in
  humans. The Journal of Neuroscience 35~(28), 10371--10385.

\bibitem[{Shadlen and Newsome(2001)}]{MNS-WTN:01}
Shadlen, M.~N., Newsome, W.~T., 2001. Neural basis of a perceptual decision in
  the parietal cortex (area {LIP}) of the rhesus monkey. Journal of
  Neurophysiology 86~(4), 1916--1936.

\bibitem[{Simen et~al.(2009)Simen, Contreras, Buck, Hu, Holmes, and
  Cohen}]{Simen2009}
Simen, P., Contreras, D., Buck, C., Hu, P., Holmes, P., Cohen, J.~D., 2009.
  Reward rate optimization in two-alternative decision making: Empirical tests
  of theoretical predictions. Journal of Experimental Psychology: Human
  Perception and Performance 35~(6), 1865.

\bibitem[{Smith(1990)}]{smith1990note}
Smith, P.~L., 1990. A note on the distribution of response times for a random
  walk with {G}aussian increments. Journal of Mathematical Psychology 34~(4),
  445--459.

\bibitem[{Smith(1995)}]{Smith1995}
Smith, P.~L., 1995. Psychophysically principled models of visual simple
  reaction time. Psychological review 102~(3), 567.

\bibitem[{Smith(2000)}]{Smith2000}
Smith, P.~L., 2000. Stochastic dynamic models of response time and accuracy: A
  foundational primer. Journal of Mathematical Psychology 44~(3), 408 -- 463.

\bibitem[{Smith and Ratcliff(2009)}]{Smith2009}
Smith, P.~L., Ratcliff, R., 2009. An integrated theory of attention and
  decision making in visual signal detection. Psychological review 116~(2),
  283.

\bibitem[{Srivastava et~al.(2016)Srivastava, Holmes, and
  Simen}]{VS-PH-PS:14-arxiv}
Srivastava, V., Holmes, P., Simen, P., 2016. Explicit moments of decision times
  for single- and double-threshold drift-diffusion processes. Journal of
  Mathematical Psychology.\; In press.

\bibitem[{Townsend and Ashby(1983)}]{Townsend1983}
Townsend, J.~T., Ashby, F.~G., 1983. Stochastic modeling of elementary
  psychological processes. Cambridge University Press.

\bibitem[{Usher and McClelland(2001)}]{Usher2001}
Usher, M., McClelland, J., 2001. The time course of perceptual choice: the
  leaky, competing accumulator model. Psychological Review 108~(3), 550.

\bibitem[{Verdonck et~al.(2015)Verdonck, Meers, and Tuerlinckx}]{Verdonck2015}
Verdonck, S., Meers, K., Tuerlinckx, F., 2015. Efficient simulation of
  diffusion-based choice rt models on {CPU} and {GPU}. Behavior Research
  Methods, 1--15.

\bibitem[{Voskuilen et~al.(2016)Voskuilen, Ratcliff, and Smith}]{Voskuilen2016}
Voskuilen, C., Ratcliff, R., Smith, P.~L., 2016. Comparing fixed and collapsing
  boundary versions of the diffusion model. Journal of Mathematical Psychology
  73, 59--79.

\bibitem[{Wagenmakers et~al.(2007)Wagenmakers, Van Der~Maas, and
  Grasman}]{Wagenmakers2007}
Wagenmakers, E.-J., Van Der~Maas, H.~L., Grasman, R.~P., 2007. An
  {EZ}-diffusion model for response time and accuracy. Psychonomic Bulletin \&
  Review 14~(1), 3--22.

\bibitem[{Wang(2002)}]{XJW:02}
Wang, X.-J., 2002. Probabilistic decision making by slow reverberation in
  cortical circuits. Neuron 36~(5), 955--968.

\bibitem[{Webb(2015)}]{Webb2015}
Webb, R., Jul. 2015. The dynamics of stochastic choice. Working paper:
  available at SSRN 2226018.

\bibitem[{White et~al.(2011)White, Ratcliff, and Starns}]{CNW-RR-JFS:11}
White, C.~N., Ratcliff, R., Starns, J.~J., 2011. Diffusion models of the
  flanker task: Dzhaniscrete versus gradual attentional selection. Cognitive
  Psychology 63~(4), 210--238.

\bibitem[{Wong and Wang(2006)}]{KFW-XJW:06}
Wong, K.-F., Wang, X.-J., 2006. A recurrent network mechanism of time
  integration in perceptual decisions. The Journal of Neuroscience 26~(4),
  1314--1328.

\bibitem[{Zhang et~al.(2014)Zhang, Lee, Vandekerckhove, Maris, and
  Wagenmakers}]{zhang2014time}
Zhang, S., Lee, M.~D., Vandekerckhove, J., Maris, G., Wagenmakers, E.-J., 2014.
  Time-varying boundaries for diffusion models of decision making and response
  time. Frontiers in Psychology 5~(1364).

\end{thebibliography}
}

\appendix

\section{Alternative expression for FPT density of {\sfr the single stage model}} \label{app:fpt-density}
An alternative expression to the FPT density~\eqref{eq:fpt-density-ddm} 
that can be obtained by solving the Fokker-Planck equation \citep{Feller1968} is: 
\begin{multline*}
f(t;x_0,\theta)
=  \frac{\pi \sigma^2}{4 \zeta^2} \exp\Big({-\frac{\mu^2 t}{2 \sigma^2}}\Big) \sum_{n=1}^{+\infty} (-1)^{n-1} n\exp\Big({- \frac{n^2 \pi^2 \sigma_1^2 t}{8 \zeta^2}}\Big)\bigg(\exp\Big({\frac{\mu(\zeta- x_0)}{\sigma^2}}\Big)\sin \Big( \frac{n \pi(\zeta+x_0) }{2 \zeta}\Big) \\
 +  \exp\Big({-\frac{\mu(\zeta+ x_0)}{\sigma^2}}\Big)\sin \Big( \frac{n \pi(\zeta-x_0) }{2 \zeta}\Big)\bigg).
\end{multline*}

\section{Derivation of expressions in \S\ref{subsec:i-ddm}}
\label{app:ddm-i}

We first establish \eqref{eq:error-rate-deadline}. 
First consider the case $\mu_i >0$. 
Let $\seqdef{\mc F_t^i}{t\ge t_{i-1}}$ be the filtration defined by the evolution of the {\sfr multistage process}~\eqref{eq:multistage-ddm} until time $t$ conditioned on $\tau > t_{i-1}$.
A filtration can be thought of as an increasing sequence of available information.  
 For some $s \in (t_{i-1}, t)$, as shown in \S\ref{sec:review},  $\seqdef{\exp({-2s_i x(t)})}{t \ge t_{i-1}}$ is a martingale, i.e., 
$\expt[\exp({-2s_i x(t)}) | \mc F_s^i] = \exp({-2s_i x(s)})$.

Furthermore, $\hat \tau_i:= \min\{\tau_i, t_i\}$ is a stopping time. Therefore, it follows from optional sampling theorem that
\begin{align*}
\expt [ \exp({-2s_i X_{i-1}})] &= \expt[\exp({ -2 s_i x(\hat \tau_i)})]\\
&= \expt[\exp({ -2 s_i x(\tau_i)})| \tau_i \le t_i] \prob [ \tau_i \le t_i ] + \expt [\exp({ -2s_i X_i})] \prob [ \tau_i > t_i ]] \\
&= \big( \exp({-2s_i \zeta}) (1- \er_i) +  \exp({2s_i\zeta}) \er_i \big) \prob [ \tau_i \le t_i ] + \expt [\exp({ -2s_i X_i})] \prob [ \tau_i > t_i ]].
\end{align*}
Solving the above equation for $\er_i$, we obtain the desired expression.

For $\mu_i =0$, , as shown in \S\ref{sec:review}, $\seqdef{x(t)}{t\ge t_i}$ is a martingale. Therefore, applying the optional sampling theorem, we obtain
\begin{align*}
\expt[X_{i-1}] &= \expt[x(\hat \tau_i)] \\
& = \expt[x(\tau_i)| \tau_i \le t_i] \prob [ \tau_i \le t_i ] + \expt[X_i] \prob [ \tau_i >t_i ]] \\
& = (1-2 \er_i) \zeta  \prob [ \tau_i \le t_i ] + \expt[X_i] \prob [ \tau_i >t_i ]. 
\end{align*}
Solving the above equation for $\er_i$, we obtain the desired expression.

The formulas \eqref{eq:fpt-density-correct} and \eqref{eq:fpt-density-incorrect} immediately follow from applying expectation to~\eqref{eq:cond-fpt-ddm-1}~and~\eqref{eq:cond-fpt-ddm-2}, respectively.

To establish \eqref{eq:decision-time-deadline}  for $\mu_i \ne 0$, we note from \S\ref{sec:review} that for the $i$-th {\sfr stage}, $\seqdef{x(t)-\mu_i t}{t\ge t_i}$ is a martingale. Therefore, applying the optional sampling theorem, we obtain
\begin{align*}
\expt[X_{i-1}] - \mu_i t_{i-1}&= \expt[x(\hat \tau_i) - \mu_i \hat \tau_i ] \\
& =  \expt[x(\tau_i) - \mu_i  \tau_i | \tau_i\le t_i] \prob [ \tau_i \le t_i ] + \expt[(X_i -\mu_i t_i ] \prob [ \tau_i >t_i ]\\
&  = \big( \zeta(1-\er_i) -\zeta \er_i - \mu_i \expt[\tau_i | \tau_i \le t_i] \big)\prob [ \tau_i \le t_i ] + (\expt[X_i]  -\mu_i t_i) \prob [ \tau_i  > t_i ]. 
\end{align*}
Solving the above equation for $\expt[\tau_i | \tau_i \le t_i]$ yields the desired expression. 

For $\mu_i=0$, we note from \S\ref{sec:review} that $\seqdef{x(t)^2 - \sigma_i^2 t}{t\ge t_{i-1}}$ is a martingale. Therefore, applying the optional sampling theorem, we obtain
\begin{align*}
\expt[X_i^2] - \sigma_i^2 t_{i-1}& =   \expt[x(\hat \tau)^2 - \sigma_i^2 \hat \tau_i] \\
& = \expt[x(\tau_i)^2 - \sigma_i^2 \tau_i|\tau_i \le t_i] \prob [ \tau_i \le t_i ] + \expt[X_i^2 - \sigma_i^2 t_i] \prob [ \tau_i >t_i ]]\\
& = (\zeta^2  -\sigma_i^2 \expt[\tau_i| \tau_i \le t_i]) \prob [ \tau_i \le t_i ] + (\expt[X_i^2] - \sigma_i^2 t_i) \prob [ \tau_i >t_i ].
\end{align*}
Solving the above equation for $\expt[\tau_i | \tau_i \le t_i]$ yields the desired expression. 

Next, we need to establish that the Laplace transform of the density for the FPT for a particular decision made before $t_i$ is
\begin{multline*}
\expt[\exp({-\alpha \tau_i}) \bs{1}(x(\tau_i)=\zeta | \tau_i \le t_i)]\\
 = \frac{\exp({-\alpha t_{i-1}})\expt[ \supscr{T}{mgf}_+(\mu_i, \sigma_i, \zeta, X_{i-1}, \alpha)] - \exp({-\alpha  t_i}) \expt[\supscr{T}{mgf}_+(\mu_i, \sigma_i, \zeta, X_i, \alpha)]\prob [  \tau _i> t_i ]]}{\prob [ \tau_i \le t_i ]}, 
\end{multline*}
\begin{multline*}
\expt[\exp({-\alpha \tau_i}) \bs{1}(x(\tau_i)=-\zeta | \tau_i \le t_i)]  \\
=\frac{ \exp({-\alpha t_{i-1}}) \expt[ \supscr{T}{mgf}_-(\mu_i, \sigma_i, \zeta, X_{i-1}, \alpha)] - \exp({-\alpha t_i}) \expt[\supscr{T}{mgf}_-( \mu_i, \sigma_i, \zeta, X_i, \alpha)\prob [  \tau_i > t_i ]]}{\prob [ \tau_i \le t_i ]}.
\end{multline*}
To establish this, we consider the stochastic process $\exp({\lambda x(t) - \lambda \mu_i t - \lambda^2 \sigma^2_i t/2})$. From \S\ref{sec:review}, we note that $\seqdef{\exp({\lambda x(t) - \lambda \mu_i t - \lambda^2 \sigma^2_i t/2})}{t\ge t_{i-1}}$ is a martingale for each $\lambda \in \real$, i.e., 
\[
\expt[\exp({\lambda x(t) - \lambda \mu_i t - \lambda^2 \sigma^2_i t/2})| \mc F_s^i] = \exp({\lambda x(s) -\lambda \mu_i s - \lambda^2 \sigma_i^2 s/2}). 
\]
We choose two particular values of $\lambda$:
\[
\lambda_1 = \frac{ -\mu_i -\sqrt{\mu_i^2 + 2 \alpha \sigma_i^2}}{\sigma_i^2}, \quad \text{and} \quad 
\lambda_2 = \frac{ -\mu_i +\sqrt{\mu_i^2 + 2 \alpha \sigma_i^2}}{\sigma_i^2}. 
\]
Note that for $\lambda \in \{\lambda_1, \lambda_2\}$, $\lambda \mu_i t + \lambda^2 \sigma^2_i t/2 = \alpha$. Therefore, stochastic processes $\seqdef{\exp({\lambda_1 x(t) - \alpha t})}{t\ge 0}$ and $\seqdef{\exp({\lambda_2 x(t) - \alpha t})}{t\ge t_{i-1}}$ are martingales. Now applying the optional sampling theorem, we obtain
\begin{multline}\label{eq:mgf-1}
\expt[\exp({\lambda_1 X_{i-1} -\alpha t_{i-1}})]  = \expt[ \exp({\lambda_1 x(\hat \tau) -\alpha \hat \tau}) ]  
=\exp({\lambda_1 \zeta}) \expt[\exp({-\alpha \tau_i}) \bs 1(x(\tau_i)=\zeta\; \& \; \tau_i \le t_i)]  \\ + 
\exp({- \lambda_1 \zeta}) \expt[\exp({-\alpha \tau}) \bs 1(x(\tau_i)=-\zeta\; \& \; \tau_i \le t_i)] +
\exp({-\alpha t_i}) \expt[\exp({\lambda_1 X_i})] \prob [ \tau_i > t_i ]. 
\end{multline}
Similarly, 
\begin{multline}\label{eq:mgf-2}
\expt[\exp({\lambda_2 X_{i-1} -\alpha t_{i-1}})] =  \exp({\lambda_2 \zeta}) \expt[\exp({-\alpha \tau_i}) \bs 1(x(\tau)=\zeta\; \& \; \tau_i \le t_i)] \\
+ \exp({- \lambda_2 \zeta}) \expt[\exp({-\alpha \tau_i}) \bs 1(x(\tau)=-\zeta\; \& \; \tau_i \le t_i)] +
\exp({-\alpha t_i}) \expt[\exp({\lambda_2 X_i})] \prob [ \tau_i > t_i ]]. 
\end{multline}
Equations~\eqref{eq:mgf-1}~and~\eqref{eq:mgf-2} are two simultaneous equations in two unknowns $\expt[\exp({-\alpha \tau_i}) \bs 1(x(\tau_i)=\zeta\; \& \; \tau_i \le t_i)]$ and $\expt[\exp({-\alpha \tau_i}) \bs 1(x(\tau_i)=-\zeta\; \& \; \tau_i \le t_i)]$. Solving for these unknowns, we obtain
\begin{multline*}
\expt[\exp({-\alpha \tau_i}) \bs 1(x(\tau_i)=\zeta\; \& \; \tau_i \le t_i)] = \frac{1}{\exp({2\lambda_1 \zeta}) - \exp({2\lambda_2 \zeta})} \bigg(
\exp({-\alpha t_{i-1}})\expt[\exp({\lambda_1(X_{i-1} +\zeta)})
\\ - \exp({\lambda_2(X_{i-1}+\zeta)})] 
 - \exp({-\alpha t_i}) \expt[\exp({\lambda_1(X_i +\zeta)}) - \exp({\lambda_2(X_i+\zeta)})] \prob [ \tau_i > t_i ]\bigg), 
\end{multline*}
and
\begin{multline*}
 \expt[\exp({-\alpha \tau_i}) \bs 1(x(\tau_i)= -\zeta\; \& \; \tau_i \le t_i)] = \frac{1}{\exp({-2\lambda_1 \zeta}) - \exp({-2\lambda_2 \zeta})}
\bigg(\exp({-\alpha t_{i-1}}) \expt[\exp({-\lambda_1(\zeta- X_{i-1})})
\\
-  \exp({- \lambda_2(\zeta- X_{i-1})})] - \exp({-\alpha t_i}) \expt[(\exp({- \lambda_1(\zeta- X_i )}) - \exp({- \lambda_2(\zeta- X_i)})) ]\prob [ \tau_i> t_i ] \bigg). 
\end{multline*}
Simplifying these expressions, we obtain the desired expression. 

Finally, \eqref{eq:decision-time-correct-deadline} and \eqref{eq:decision-time-incorrect-deadline} follow from differentiating the Laplace transform with respect to $-\alpha$, and then evaluating at $\alpha =0$.

\section{Performance metrics for the overall {\sfr multistage} process}
\label{app:mult-ddm}

We start by establishing~\eqref{eq:multistage-fpt}. Since $t \in (t_{k-1}, t_k]$, 
\begin{align*}
\prob [ \tau \le t ] & = \prob [ \tau \le t  \; \&\;  \tau \le  t_{k-1} ] + \prob [ \tau \le t  \; \&\;  \tau > t_{k-1} ] \\
& = \prob [ \tau \le t_{k-1} ] + \prob [ \tau \le t | \tau > t_{k-1} ] \prob [ \tau > t_{k-1} ] \\
& = 1 - \prod_{i=1}^{k-1} \prob [ \tau > t_i | \tau > t_{i-1} ] +  \prob [ \tau_k \le t ] \prod_{i=1}^{k-1} \prob [ \tau > t_i | \tau > t_{i-1} ] \\
& = 1 - \prod_{i=1}^{k-1} \prob [ \tau_i > t_i ] +  \prob [ \tau_k \le t ] \prod_{i=1}^{k-1} \prob [ \tau_i > t_i ]. 
\end{align*}

We now establish~\eqref{eq:multistage-DT}. We note that
\begin{align*}
\expt[\tau] &= \sum_{i=1}^n \expt[\tau \bs 1(t_{i-1} < \tau \le  t_i )] \\
& = \sum_{i=1}^n \expt[\tau \bs 1(\tau \le  t_i ) | \tau > t_{i-1}] \prob [ \tau > t_{i-1} ] \\ 
& =  \sum_{i=1}^n \Big( \expt[\tau_i | \tau_i \le  t_i ]  \prob [ \tau_i \le  t_i  ] \prod_{j=1}^{i-1} \prob [ \tau_j > t_j ] \Big). 
\end{align*}

To establish~\eqref{eq:multistage-ER}, we note that
\begin{align*}
\er &= \sum_{i=1}^{n+1} \prob (x(\tau) = -\zeta \text{ and } t_{i-1} \le \tau < t_i ) \\
&= \sum_{i=1}^{n+1} \prob (x(\tau)  -\zeta \text{ and } \tau < t_i | \tau > t_{i-1}) \prob [ \tau > t_{i-1} ] \\
& = \sum_{i=1}^n \Big( \er_i \prob [ \tau_i < t_i  ] \prod_{j=1}^{i-1} \prob [ \tau_j > t_j ] \Big). 
\end{align*}

Equations~\eqref{eq:multistage-DT-correct}~and~\eqref{eq:multistage-DT-incorrect} follow similarly to~\eqref{eq:multistage-DT}, and 
Equations~\eqref{eq:multistage-fpt-correct}~and~\eqref{eq:multistage-fpt-incorrect} follow similarly to~\eqref{eq:multistage-fpt}.

\section{Time-varying Ornstein-Uhlenbeck model}
\label{app:ou}

In this section, we discuss how the ideas presented in \S\ref{sec:performance-metrics}~and~\S\ref{sec:thresholds} can be used to computed FPT properties for a time varying 
Ornstein-Uhlenbeck (O-U) model.

{\sfr The} O-U model captures decision-making through the first passage of trajectories of an O-U process~\citep{DRC-HDM:65} through two thresholds. 
Our calculations for the {\sfr multistage process} also help analyze the Ornstein-Uhlenbeck (O-U) model.
Similar to the {\sfr multistage process}, the $n$-stage O-U process with piecewise constant parameters is defined by
\begin{equation}\label{eq:multistage-o-u}
d x(t) = \mu(t) d t - \lambda(t) x(t) dt + \sigma(t) d W(t), \quad x(t_0) = x_0, 
\end{equation}
where 
\begin{align*}
\mu(t) &= \mu_i,  \quad \text{for } t_{i-1} \le t < t_i, \\
\sigma(t) &=   \sigma_i, \quad  \text{for } t_{i-1} \le t < t_i, \\
\lambda(t) &=   \lambda_i,  \quad  \text{for } t_{i-1} \le t < t_i,
\end{align*}
for each $i \in \until{n}$. Due to the extra $- \lambda(t) x(t) dt $ term, the O-U process is a leaky integrator, while the {\sfr original single stage process} is a perfect integrator ($\lambda (t) = 0$). 

Here, leaky integration means that as the noisy signal is integrated in time with exponentially increasing ($\lambda<0$) or decreasing ($\lambda>0)$ weights on past observations. Such exponential weights lead to `recency' or `decay' effects, i.e.,  the earlier stages (or late stages) may have greater influence on the ultimate decision, whereas with the {\sfr single stage process}, all of the signal throughout the entire decision period is weighed equally.

\subsection{The O-U process as a transformation of the Wiener process} \label{subsec:trans-weiner-o-u}

In this section we show how our calculations for the {\sfr multistage process} can be easily applied to decision models driven by O-U processes via a transformed Wiener process~\cite[\S5.9]{DRC-HDM:65}.

Consider the single-stage O-U process
\begin{equation}\label{eq:o-u-process}
d x(t) = \mu dt - \lambda x(t) dt + \sigma d W(t), \; x(t_0) =x_0.
\end{equation}

The O-U process~\eqref{eq:o-u-process} can be written as a
time-varying location-scale transformation of the Wiener process~\cite[\S5.9]{DRC-HDM:65}, i.e., 
\begin{equation}\label{eq:ou-transform}
x(t) = \frac{\mu}{\lambda}(1- \exp({-\lambda t})) + \exp({-\lambda t}) x_0 + \exp({-\lambda t}) W \Big( \frac{\sigma^2 (\exp({2\lambda t}) -1)}{2 \lambda} \Big).
\end{equation}

In order to derive~\eqref{eq:ou-transform}, note that for O-U process~\eqref{eq:o-u-process} 
\[
x(t) = x(0) \exp(-\lambda t) + \frac{\mu}{\lambda}  (1- \exp(-\lambda t) ) + \sigma \int_0^t \exp(-\lambda s) dW(s),
\]
Note that the stochastic process $\sigma \int_0^t \exp(-\lambda s) dW(s)$ is equivalent to the stochastic process $\sigma W((1-\exp(-2\lambda t))/2\lambda)$ in the sense of distribution. Furthermore, $\sigma W((1-\exp(-2\lambda t))/2\lambda)$ is equivalent to the stochastic process $e^{-\lambda t} W(\sigma^2 (\exp(2\lambda t) -1 )/2\lambda)$ in a similar sense. This means that for each realization of the process $\sigma \int_0^t \exp(-\lambda s) dW(s)$ there exists an identical realization of the process $e^{-\lambda t} W(\sigma^2 (\exp(2\lambda t) -1 )/2\lambda)$.

If we define the transformed time by $u(t) : = \frac{\sigma^2 (\exp({2\lambda t}) -1)}{2 \lambda}$ so that $t = \frac{1}{2 \lambda_1} \log (1 + \frac{2 \lambda u}{\sigma^2})$. Then  
\[
x_0 + W(u(t))  = \Big(x(t) - \frac{\mu}{\lambda} \Big) \exp({\lambda t}) +  \frac{\mu}{\lambda}  = \Big(x(t) - \frac{\mu}{\lambda} \Big) \sqrt{1 + \frac{2 \lambda u(t)}{\sigma^2}} +  \frac{\mu}{\lambda}.
\]
We refer to this process as a Wiener process evolving on exponential time scale.

\subsection{First passage time of the O-U process}
We now consider the first passage time of the O-U process~\eqref{eq:o-u-process} with respect to symmetric thresholds $\pm \zeta$. 
If $x(t) = \pm \zeta$, we have $x_0 + W(u) = \big(\pm \zeta -
\frac{\mu}{\lambda} \big) \sqrt{1 + \frac{2 \lambda
    u}{\sigma^2}} + \frac{\mu}{\lambda}$. We denote this last
quantity by $ \zeta^{\pm}(u)$. Consequently, the FPT for $x(t)$ with
respect to thresholds $\pm \zeta$ is a continuous transformation of the
FPT of a Wiener process starting at $x_0$ and evolving on the
transformed time $u$ with respect to time-varying thresholds at
$\zeta^\pm(u)$.  Since $u$ is a monotonically increasing function,
the distribution of the first passage time $\tau$ can be obtained
from $u(\tau)$, the FPT distribution of the Wiener
process. Furthermore, the lower threshold hitting probabilities of the
two processes are the same.

Note that in transforming the FPT problem for the O-U process~\eqref{eq:o-u-process} to the FPT problem for the Wiener process evolving on exponential time scale, removes all the parameters from the underlying process $x_0+W(u)$ and puts them in thresholds $\zeta^\pm(u)$ and exponential time scale $u$. In addition to {\sfr utilizing the results of \S \ref{sec:performance-metrics}} to multistage O-U processes, the above transformation is also helpful is speeding up Monte-Carlo simulations of the O-U process. Since the transformed process evolves on exponential time scale, the Monte-Carlo simulations with transformed process should roughly take time that is a logarithmic function of time taken by the O-U process~\eqref{eq:o-u-process}.

Computation of FPT distributions for the Wiener process with
time-varying thresholds is, to our knowledge, not analytically
tractable. However, the time-varying thresholds can be approximated by
piecewise constant time-varying thresholds and approximate FPT
distributions can be computed using the {\sfr multistage model}. While the thresholds
$\zeta^\pm$ are asymmetric for the transformed process (i.e. $\zeta^+
\ne -\zeta^-$), unlike in the case described for the {\sfr multistage process}; such a
case can be easily handled by replacing the expression
in~\eqref{eq:fpt-density-ddm}~and~\eqref{eq:dist-deadline} with
corresponding expressions for asymmetric thresholds (see
\citealp{Douady1999, Borodin2002}).

\subsection{Approximate computation of the FPT distribution of the multistage O-U process}

Similar to the transformation described in \ref{subsec:trans-weiner-o-u}, the multistage O-U process~\eqref{eq:multistage-o-u} for $t \in [t_{i-1}, t_i)$ can be written as
\begin{multline}\label{eq:i-th-o-u-solution}
x(t) = \frac{\mu_i}{\lambda_i}(1- \exp({-\lambda_i (t- t_{i-1})})) + \exp({-\lambda_i (t- t_{i-1})}) x(t_{i-1}) \\
+ \exp({-\lambda_i (t- t_{i-1})}) W \Big( \frac{\sigma_i^2 (\exp({2\lambda_i (t- t_{i-1})}) -1)}{2 \lambda_i} \Big).
\end{multline}
Let $u_i(t) = {\sigma_i^2 (\exp({2\lambda_i (t- t_{i-1})}) -1)}/{2 \lambda_i}$.
Also, let $\tau_i$ and $X_{i-1}$, $i \in \until{n}$ be defined similarly to the multistage model.
 Then, conditioned on a realization of $X_{i-1}$, the FPT problem of the $i$-th stage O-U process can be equivalently written as the FPT problem of the Wiener process $X_{i-1} + W(u_i(t))$ with respect to thresholds
\[
\zeta_i^\pm (u_i(t)) = \Big(\pm \zeta -\frac{\mu_i}{\lambda_i} \Big) \sqrt{1+ \frac{2 \lambda_i u_i(t)}{\sigma_i^2} } + \frac{\mu_i}{\lambda_i}.
\]
We can approximate each stage of the O-U process~\eqref{eq:multistage-o-u} by {\sfr a multistage process} representing the above Wiener process with time varying thresholds. {\sfr This sequence of multistage processes is itself a larger multistage process} that approximates~\eqref{eq:multistage-o-u} and its FPT distribution can be computed using the method developed in \S\ref{sec:performance-metrics}. Note that this method only yields FPT distributions. The expected decision times and probability of hitting a particular threshold can be computed by integrating these distributions. 

\end{document}